\documentclass{amsart}
\usepackage{amsthm}
\usepackage[dvips]{graphicx}

\makeatletter

\newcommand{\bmath}[1]{\mbox{\mathversion{bold}$#1$}}
\newcommand{\smath}[1]{\mbox{\mathversion{bold}\scriptsize{$#1$}}}
\newcommand{\C}{\bmath{C}}
\newcommand{\subC}{\smath{C}}
\newcommand{\Z}{\bmath{Z}}
\newcommand{\R}{\bmath{R}}

\newcommand{\CP}{\bmath{C\!P}}
\newcommand{\SL}{\operatorname{SL}}
\newcommand{\SU}{\operatorname{SU}}
\newcommand{\PSL}{\operatorname{PSL}}
\newcommand{\PSU}{\operatorname{PSU}}
\newcommand{\cmcone}{\mbox{\rm CMC}-$1$}
\newcommand{\TA}{\operatorname{TA}}
\newcommand{\trans}{{}^t_{}\!}
\newcommand{\gO}{\mathbf{O}}
\newcommand{\gI}{\mathbf{I}}

\newcommand{\gGamma}{\bmath{\Gamma}}
\newcommand{\ord}{\operatorname{ord}}
\newcommand{\id}{\operatorname{identity}}
\newcommand{\Hyp}{\mathcal{H}}
\renewcommand{\Re}{\operatorname{Re}}
\newcommand{\Per}{\operatorname{Per}}
\newcommand{\Res}{\operatornamewithlimits{Res}}
\renewcommand{\Im}{\operatorname{Im}}
   \newtheorem{theorem}{Theorem}[section]
   \newtheorem{proposition}[theorem]{Proposition}
   
   \newtheorem{lemma}[theorem]{Lemma}
   \newtheorem*{fact}{Fact}
 \theoremstyle{definition}
   \newtheorem{defn}{Definition}

 \theoremstyle{remark}
   \newtheorem{example}[theorem]{Example}

\numberwithin{equation}{section}

\def\subsubsection{\@startsection{subsubsection}{3}%
  {\parindent}{.5\linespacing\@plus.7\linespacing}{-.5em}%
  {\bf}}
\def\paragraph{\@startsection{paragraph}{4}%
  {\z@}{\z@}{-\fontdimen2\font}%
  \normalfont\itshape}
\title[CMC-1 surfaces of low total curvature I]{
   Mean curvature 1 surfaces in hyperbolic
   3-space with low total curvature I
}
\dedicatory{Dedicated to Katsuhiro Shiohama on the
occasion of his sixtieth birthday.}
\date{July 12, 2000}
\author{Wayne Rossman}
\author{Masaaki Umehara}
\author{Kotaro Yamada}
\address[Rossman]{%
   Department of Mathematics, Faculty of Science,
   Kobe University,
   Rokko, Kobe 657-8501, Japan%
}
\email{wayne@math.kobe-u.ac.jp}
\address[Umehara]{%
   Department of Mathematics, Faculty of Science,
   Hiroshima University,
   Higashi-Hiroshima 739-8526, Japan%
}
\email{umehara@math.sci.hiroshima-u.ac.jp}
\address[Yamada]{%
   Faculty of Mathematics,
   Kyushu University 36, 6-10-1
   Hakozaki, Higashi-ku, Fukuoka 812-8185, Japan%
}
\email{kotaro@math.kyushu-u.ac.jp}
\begin{document}
\begin{abstract}
A complete surface of constant mean curvature $1$ (CMC-$1$) in 
hyperbolic $3$-space with constant curvature $-1$
has two natural notions of ``total curvature''---
one is the {\em total absolute curvature} which is the integral over 
the surface of the absolute value of the Gaussian curvature, 
and the other is the {\em dual total absolute curvature} which is the 
total absolute curvature of the dual CMC-$1$ surface.  In this 
paper, we completely classify CMC-$1$ surfaces with dual total 
absolute curvature at most $4\pi$.
Moreover, we give new examples  and partially classify CMC-$1$ 
surfaces with dual total absolute curvature at most $8\pi$.
\end{abstract}
\maketitle
With the developments of the last decade on constant mean curvature $1$
(\cmcone) surfaces in hyperbolic $3$-space $H^3$ (the complete
simply-connected $3$-manifold of constant sectional curvature $-1$), and
with so many examples now known, it is a natural next step to classify
all such surfaces with low total absolute curvature.  

As \cmcone{} surfaces in $H^3$ share quite similar properties with 
minimal surfaces in Euclidean $3$-space $\R^3$, let us first comment
that the total absolute curvature of a minimal surface in $\R^3$ is
equal to the area (counted with multiplicity) of the Gauss image of the
surface, and that complete minimal surfaces in $\R^3$ with total 
curvature at most $8\pi$ have been classified. 
(See Lopez~\cite{Lopez} and also Table~\ref{tab:minimal}.)  
Furthermore, as the Gauss map of a complete conformally parametrized 
minimal surface is holomorphic, and has a well-defined limit at each end
when the surface has finite total curvature, 
the area of the Gauss image must be an integer multiple of $4\pi$. 

The question of classifying low total curvature \cmcone{} surfaces in
$H^3$ is analogous --- however, unlike minimal surfaces in $\R^3$,
\cmcone{} surfaces in $H^3$ have two Gauss maps: 
the hyperbolic Gauss map $G$ and the secondary Gauss map $g$.  
So there are two ways to pose the question in $H^3$, with two very
different answers.
One way is to consider the true total absolute curvature, 
which is the area of the image of $g$, but since $g$ might not be 
single-valued on the surface, the total curvature might not be an
integer multiple of $4\pi$.  
This allows for many more possibilities and makes the problem more
difficult than for minimal surfaces in $\R^3$.  
The authors take up this question in a separate paper \cite{ruy4}.  

The second way, which is the theme of this paper, is to study the area
of the image of $G$, 
which we call the {\it dual\/} total absolute curvature, as it is the
true total curvature of the dual \cmcone{} surface (which we define in 
Section~\ref{sec:summary}) in $H^3$.  
This way has the advantage that $G$ is single-valued on the surface, and so
the dual total curvature is always an integer multiple of $4\pi$, 
like the case of minimal surfaces in $\R^3$.  
Furthermore, the dual total curvature satisfies not only the Cohn-Vossen
inequality, but also the hyperbolic analogue of the Osserman inequality 
(which cannot be said about the true total curvature) 
\cite{uy5,Yu2} 
(see also \eqref{eq:osserman} in Section~\ref{sec:prelim}).  
So the dual total curvature shares more properties with the total
curvature of minimal surfaces in $\R^3$, motivating our interest in it.  

In this paper, we classify \cmcone{} surfaces with 
dual total absolute curvature at most $4\pi$, and we go much of the way 
toward classifying \cmcone{} surfaces with 
dual total absolute curvature at most $8\pi$ (as a first step to a 
full classification of the $8\pi$ case).  
In Section~\ref{sec:summary}, we give a summary of the results, and 
in Section~\ref{sec:prelim} we give preliminaries for the latter sections.
The classification of \cmcone{} surfaces with dual total absolute
curvature less than or equal to $4\pi$ is given in Section~\ref{sec:4-pi}.
Surfaces with dual total absolute curvature $8\pi$ are discussed in
Section~\ref{sec:8-pi} --- and there we find new examples, we classify 
certain cases, and we show nonexistence in certain other cases.  
In Section~\ref{sec:deform}, from deformations of corresponding minimal 
surfaces in $\R^3$, we produce two classes of new \cmcone{} surfaces with 
dual total absolute curvature $8\pi$.  
For the readers' convenience, we attach 
Appendix~\ref{app:log} to explain the computation of log-term
coefficients of second order linear ordinary differential equations with
regular singularities.
\section{Summary of the results}
\label{sec:summary}

To state our results precisely, we begin with some notations.  
Let $f\colon{}M\to H^3$ be a complete conformal \cmcone{} immersion of 
a Riemann surface $M$ into $H^3$.  
By Bryant's representation formula, there is a holomorphic null
immersion $F\colon{}{\widetilde M}\to\SL(2,\C)$ such that 
$f=FF^{*}$, where $\widetilde M$ is the universal cover of $M$ and
$F^{*}=\trans\overline{F}$.  (``null'' means $\det (F^{-1}dF)=0$.)  
Here, we consider $H^3=\SL(2,\C)/\SU(2)=\{aa^{*}\,|\,a\in\SL(2,\C)\}$
\cite{Bryant,uy1}.
We call $F$ the {\em lift\/} of $f$, and $F$ satisfies 
\begin{equation} \label{eq:F}
  dF = F 
  \begin{pmatrix} g &  -g^2 \\ 1 & -g\hphantom{^2} \end{pmatrix} 
  \frac{Q}{dg}
\end{equation} 
on $\widetilde M$, where $g$ (the {\it secondary Gauss map\/}) 
is a meromorphic function defined on $\widetilde M$ and $Q$ (the 
{\it Hopf differential\/}) is a holomorphic $2$-differential on $M$.
Then the induced metric $ds^2$ and complexification of the 
second fundamental form $h$ are 
\[
    ds^2 = (1+|g|^2)^2 \left| \frac{Q}{dg} \right|^2 \; , \qquad
    h= -Q-\overline{Q}+ds^2 \; . 
\]  
By \eqref{eq:F}, the secondary Gauss map satisfies
\[
    g = -\frac{dF_{12}}{dF_{11}} = -\frac{dF_{22}}{dF_{21}} \; , 
    \quad \text{where} \quad
    F(z) = \begin{pmatrix}F_{11}(z) & F_{12}(z) \\ 
                          F_{21}(z) & F_{22}(z) \end{pmatrix}\; . 
\]
The map $g$ is determined uniquely up to a M\"obius transformation 
$g\mapsto a\star g$ by $a\in\SU(2)$, where, for general 
$a=(a_{ij})\in\SL(2,\C)$, we denote
\[
   a\star g :=\frac{a_{11}g + a_{12}}{a_{21}g+ a_{22}}\;.
\]
The {\it hyperbolic Gauss map\/} $G$ of $f$ is defined by
\[
   G = \frac{dF_{11}}{dF_{21}} = \frac{dF_{12}}{dF_{22}}\;,
\]
which can be interpreted as stereographic projection of the endpoints in
the sphere at infinity of $H^3$ of the oriented normal geodesics
emanating from the surface.
In particular, $G$ is a meromorphic function on $M$.  

The inverse matrix $F^{-1}$ is also a holomorphic null immersion,
and produces a new \cmcone{} immersion 
$f^{\#}=F^{-1}(F^{-1}{})^{*}\colon{}\widetilde M\to H^3$, 
called the {\em dual\/} of $f$ \cite{uy5}.  
The induced metric $ds^2{}^{\#}$ and the Hopf differential $Q^{\#}$ of 
$f^{\#}$ are 
\begin{equation}\label{eq:fund-forms}
    ds^2{}^{\#}= (1+|G|^2)^2 \left|\frac{Q}{dG}\right|^2\;,
    \qquad
    Q^{\#}= - Q \; . 
\end{equation}
So $ds^2{}^{\#}$ and $Q^{\#}$ are well-defined on $M$ itself, 
even though $f^{\#}$ might be defined only on $\widetilde M$.
This duality between $f$ and $f^{\#}$ interchanges the roles of the 
hyperbolic Gauss map $G$ and secondary Gauss map $g$.  In particular, one has 
\begin{equation} \label{eq:Fsharp}
  dF\, \cdot \, F^{-1}= -(F^{-1})^{-1}d(F^{-1})=
    \begin{pmatrix} G & -G^2 \\ 1 & -G\hphantom{^2} \end{pmatrix}
     \frac{Q}{dG}\;.
\end{equation} 
Hence $dFF^{-1}$ is single-valued on $M$, whereas $F^{-1}dF$ generally is not.

Since  $ds^2{}^{\#}$ is single-valued on $M$,
we can define the {\em dual total absolute curvature\/} 
\[
    \TA(f^{\#}) := \int_M (-K^{\#})\,dA^{\#},
\]
where $K^{\#}$ ($\leq 0$) and $dA^{\#}$ are the Gaussian curvature and
area element of $ds^2{}^{\#}$, respectively.  As 
\begin{equation}\label{eq:pseudo}
   d\sigma^2{}^{\#}:=(-K^{\#})ds^2{}^{\#} = 
        \frac{4\,dG\,d\overline{G}}{(1+|G|^2)^2}
\end{equation}
is a pseudo-metric of constant curvature $1$ with developing map $G$,
$\TA(f^{\#})$ is the area of the image of $G$ on $\CP^1=S^2$.
The following assertion is important for us:  
\begin{lemma}[\cite{uy5,Yu}] \label{lem:complete} 
  The Riemannian metric $ds^2{}^{\#}$ is complete
  $($resp. nondegen\-erate$)$
  if and only if $ds^2$ is complete $($resp.~nondegenerate$)$.  
\end{lemma}
So from now on, we suppose $f$ is complete and has $\TA(f^{\#})<+\infty$. 
By Lemma~\ref{lem:complete}, the conformal metric $ds^{2}{}^{\#}$ is
complete.  
As $\TA(f^{\#})<+\infty$, $M$ is biholomorphic  to a compact Riemann
surface $\overline{M}_\gamma$ of some genus $\gamma$ with finitely many
points excluded \cite[Theorem 9.1]{Osserman}: 
\begin{equation}\label{eq:compactify}
   M = \overline{M}_\gamma \setminus\{p_1,\dots,p_n\}\qquad 
       (p_1,\dots, p_n\in \overline{M}_\gamma) \; . 
\end{equation}
The points $p_j$ are called the {\it ends\/} of the immersion $f$.

If $G$ has an essential singularity at any end $p_j$, then
$\TA(f^{\#})=+\infty$, since $\TA(f^{\#})$ is the area of $G(M)$ in
$\CP^1=S^2$.  
Since we have assumed $\TA(f^{\#}) < +\infty$, $G$ is meromorphic on all
of $\overline M_{\gamma}$.  
In particular, $\TA(f^{\#})=4\pi\deg G\in 4\pi\Z$.  

\begin{table}
\begin{small}
  \begin{tabular}{|l||c|l|l|l|l|}
    \hline
      Type & {\footnotesize$\TA(f^{\#})$} & Reducibility & Status & 
                                \multicolumn{1}{|c|}{c.f.} \\
    \hline\hline
      $\gO(0)$  & 0 & $\Hyp^3$-reducible
                & classified$^0$
                & Horosphere    \\
    \hline
    \hline
      $\gO(-4)$ & $4\pi$ & $\Hyp^3$-reducible
                & classified
                & \parbox[t]{4.4cm}{
                    Duals of Enneper cousins \\
                    \hspace*{\fill}\cite[Example~5.4]{ruy1}}    \\
    \hline
      $\gO(-2,-2)$ & $4\pi$
                & 
               reducible
                                & classified
                & \parbox[t]{4.4cm}{Catenoid cousins and warped
                catenoid cousins with embedded ends \\
                \hspace*{\fill}
                \cite[Example~2]{Bryant},\cite{uy1},\cite{ruy4}}   \\
    \hline
    \hline
      $\gO(-5)$ & $8\pi$ & $\Hyp^3$-reducible
                & classified
                & Theorem~\ref{thm:go5/6} \\
    \hline
      $\gO(-6)$ & $8\pi$ & $\Hyp^3$-reducible
                & classified
                & Theorem~\ref{thm:go5/6} \\
    \hline
      $\gO(-2,-2)$ 
                & $8\pi$ & reducible
                & classified
                & \parbox[t]{4.4cm}{Double covers of \\
                    catenoid cousins and warped 
                    catenoid cousins with  $m=2$ in \\
                    \hspace*{\fill}\cite[Theorem~6.2]{uy1},\cite{ruy4}}
                    \\
    \hline
      $\gO(-1,-4)$ 
                & $8\pi$ & $\Hyp^3$-reducible
                & classified$^0$
                & Theorem~\ref{thm:go14} \\
    \hline
      $\gO(-2,-3)$ 
                & $8\pi$ & $\Hyp^1$-reducible
                & classified
                & 
                Theorems~\ref{thm:go23-0}, \ref{thm:go23-1} \\
    \hline
      $\gO(-2,-4)$ 
                & $8\pi$ & $\Hyp^1$-reducible
                & classified &
              Theorem~\ref{thm:go24} \\
                &  & $\Hyp^3$-reducible
                & classified &
              Theorem~\ref{thm:go24b} \\
    \hline
      $\gO(-3,-3)$ 
                & $8\pi$ & reducible
                & existence
                & Proposition~\ref{thm:go33} \\
    \hline
      $\gO(-1,-1,-2)$ 
                & $8\pi$ &  $\Hyp^3$-reducible
                & classified$^0$
                & Theorem~\ref{thm:go112}\\
    \hline
      $\gO(-1,-2,-2)$ 
                & $8\pi$ & $\Hyp^1$-reducible
                & classified
                & Theorem~\ref{thm:go122} \\

                &  & $\Hyp^3$-reducible
                & classified    
                & Theorem~\ref{thm:go122-h3} 
\\
    \hline
      $\gO(-2,-2,-2)$ 
                & $8\pi$ & irreducible
                & classified
                &  \cite[Theorem~2.6]{uy6}            \\
                &  & $\Hyp^1$-reducible
                & existence$^+$
                & Example~\ref{ex:o222red1}    \\
                &  & $\Hyp^3$-reducible
                & existence$^+$
                & Example~\ref{ex:o222red3}    \\
    \hline
      $\gI(-3)$ 
                & $8\pi$ &
                & unknown
                & \\ 
    \hline 
      $\gI(-4)$ 
                & $8\pi$ &
                & existence
                & Proposition~\ref{thm:gi4}    \\
    \hline
      $\gI(-1,-1)$ 
                & $8\pi$ &
                & unknown$^{+}$
                & Proposition~\ref{prop:data-I}    \\
    \hline
      $\gI(-2,-2)$ 
                & $8\pi$ &
                & existence
                & Genus $1$ catenoid cousins~\cite{rs}    \\
    \hline
  \end{tabular}
\end{small}
\caption{\cmcone{} surfaces in $H^3$ with 
         $\TA(f^{\#})\leq 8\pi$.  (The corresponding results for 
         minimal surfaces in $\R^3$ are shown in Table~\ref{tab:minimal}.)}
\label{tab:class}
\end{table}

Since the dual immersion has finite total curvature, the Hopf 
differential $Q^{\#}=-Q$ can be extended to $\overline M_{\gamma}$ as
a meromorphic $2$-differential~\cite[Proposition 5]{Bryant}.
Let 
\[
   d_j = \ord_{p_j}Q = \text{ order of $Q$ at the end $p_j$} 
\]
for each  $j=1,\dots,n$.  
We say that $f$ is a surface of {\em type~$\gGamma(d_1,\dots,d_n)$} if 
$M = \overline{M}_{\gamma}\setminus\{p_1,\dots,p_n\}$ and $Q$ has order 
$d_j$ at each end $p_j$.  
We use $\gGamma$ because it is the capitalized letter corresponding to 
$\gamma$, the genus of $\overline M_{\gamma}$.  
For instance, the class $\gI(-4)$ (resp.~$\gO(-2,-3)$) means the class of 
surfaces of genus $1$ (resp.~genus $0$) with $1$ end (resp.~$2$ ends) so that 
$Q$ has a pole of order $4$ at the single end (resp.~a pole of order $2$ at 
one end and order $3$ at the other).  
Then our results are shown in Table~\ref{tab:class}.
In the table, 
\begin{itemize}
  \item {\em classified\/} means the complete list of the surfaces
        in such a class is known (and this means not only that we know 
        all the possibilities for the form of the data $(G,Q)$, 
        but that we also know exactly for which $(G,Q)$ the period problems 
        of the immersions are solved).  
  \item {\em classified\/$^0$} means there exists a unique surface
        (up to isometries of $H^3$ and deformations that come from its
           reducibility).  
  \item {\em existence\/} means that examples exist, but they are not yet 
                  classified.  
  \item {\em existence{$^+$}} means that all possibilities for the 
        data $(G,Q)$ are determined in this paper,
        but the period problems are solved only for special cases.
  \item {\em unknown\/} means that neither existence nor non-existence is
        known yet.
  \item {\em unknown{$^+$}} means that all possibilities for the 
        data $(G,Q)$ are determined in this paper, but the period
	problems are still unsolved.
 \end{itemize}
Any class and type of reducibility not listed in Table \ref{tab:class} cannot 
contain surfaces with $\TA(f^{\#}) \leq 8\pi$.  For example, any irreducible 
or $\Hyp^3$-reducible surface of type $\gO(-2,-3)$ must have dual total
absolute curvature at least $12 \pi$.  
(See Section~\ref{sec:prelim} for the definitions of irreducibility, 
$\Hyp^1$-reducibility, and $\Hyp^3$-reducibility.)
\begin{table}[ht]
  \begin{tabular}{|l||c|l|l|}
    \hline
      Type & $\TA$ & The surface & \multicolumn{1}{|c|}{c.f.} \\
    \hline\hline
      $\gO(0)$  & 0 
                & Plane
                &     \\
    \hline
      $\gO(-4)$ & $4\pi$ 
                & Enneper's surface
                &    \\
    \hline
      $\gO(-5)$ & $8\pi$ 
                & 
                & \cite[Theorem~6]{Lopez}   \\
    \hline
      $\gO(-6)$ & $8\pi$ 
                & 
                & \cite[Theorem~6]{Lopez}   \\
    \hline
      $\gO(-2,-2)$ 
                & $4\pi$
                & Catenoid
                & \\
                & $8\pi$
                & Double cover of the catenoid
                & \\
    \hline
      $\gO(-1,-3)$ 
                & $8\pi$
                & 
                & \cite[Theorem~5]{Lopez}\\
    \hline
      $\gO(-2,-3)$ 
                & $8\pi$
                & 
                & \cite[Theorem~4, 5]{Lopez}\\
    \hline
      $\gO(-2,-4)$ 
                & $8\pi$
                & 
                & \cite[Theorem~5]{Lopez}\\
    \hline
      $\gO(-3,-3)$ 
                & $8\pi$
                & 
                & \cite[Theorem~4]{Lopez}\\
    \hline
      $\gO(-1,-2,-2)$ 
                & $8\pi$
                & 
                & \cite[Theorem~5]{Lopez}\\
    \hline
      $\gO(-2,-2,-2)$ 
                & $8\pi$
                & 
                & \cite[Theorem~5]{Lopez}\\
    \hline
      $\gI(-4)$ 
                & $8\pi$
                & Chen-Gackstatter surface
                & \cite[Theorem~5]{Lopez}, \cite{cg}\\
    \hline
  \end{tabular}
\vspace{0.1in}
\caption{The classification of complete minimal surfaces in $\R^3$ with 
$\TA\leq 8\pi$ (\cite{Lopez}), for comparison with Table~\ref{tab:class}.}
\label{tab:minimal}
\end{table}
\section{Preliminaries}
\label{sec:prelim}

Before we begin proving the results, we prepare some
fundamental properties and tools, which will
play important roles in the latter sections.

\subsection*{Analogue of the Osserman inequality.} 
The second and third authors showed \cite{uy5}: 
\begin{equation}\label{eq:osserman}
   \frac{1}{2\pi}\TA(f^{\#})\geq -\chi(M)+n = 2(\gamma+n-1) \; . 
\end{equation}   
Moreover, equality holds exactly when all the ends are embedded: This follows 
by noting that equality is equivalent to all ends being regular and 
embedded (\cite{uy5}), and that any embedded end must be regular (proved 
recently by Collin, Hauswirth and Rosenberg \cite{chr}).  

\subsection*{Formulas for $\mathbf{TA(\bmath{f}^{\#})}$.}  
Let $\mu_j^{\#}\in\Z$ be the branching order of $G$ at the end $p_j$
for each $j=1,\dots,n$.
Since $G$ is a ($\mu_j^{\#}+1$)-to-$1$ mapping in a neighborhood of $p_j$, 
\begin{equation}\label{eq:mu-ineq}
   \mu_j^{\#}\leq \deg G-1=\frac{1}{4\pi}\TA(f^{\#}) -1\;.
\end{equation}

The umbilic points of $f$ are the zeroes of $Q=-Q^{\#}$, which are also
the umbilic points of $f^{\#}$. 
Moreover, the order of $Q$ equals the branching order of $G$ at each point 
in $M$, since $ds^2{}^{\#}$ in \eqref{eq:fund-forms} is non-degenerate.  
Let $q_1, \dots, q_k$ be the umbilic points of $f$ and set
\begin{equation}\label{eq:QandG}
   \xi_l:= \ord_{q_l}Q =
     \left[\text{the branching order of $G$ at $q_l$}\right]
   \qquad (l=1,\dots,k)\;.
\end{equation}
The pseudometric $d\sigma^2{}^{\#}$ in \eqref{eq:pseudo} is said to have 
order $\beta$ at $p$ if it is asymptotic to 
$|z-z(p)|^{2\beta}dz d\bar z$, where $z$
is a complex coordinate around $p$. 
Then the branching order of $G$ is equal to the order of the metric
$d\sigma^2{}^{\#}$ in \eqref{eq:pseudo}, the Gauss-Bonnet theorem
implies that 
\begin{equation}\label{eq:gauss-bonnet}
   \frac{1}{2\pi}\TA(f^{\#}) = \chi(\overline M_{\gamma})+
              \sum_{j=1}^{n}\mu_j^{\#} + \sum_{l=1}^k \xi_l\;,
\end{equation}
where $\chi(\,\cdot\,)$ is the Euler characteristic.  
(This also follows from the Riemann-Hurwitz formula, 
since $\overline{M}_\gamma$ is a branched cover of $S^2$ via the map $G$.)  

Since $Q$ is a meromorphic $2$-differential, the total order of $Q$ 
satisfies 
\begin{equation}\label{eq:r-r}
   \sum_{l=1}^k\xi_l + \sum_{j=1}^n d_j = -2\chi(\overline M_{\gamma})\;.
\end{equation}
By \eqref{eq:gauss-bonnet} and \eqref{eq:r-r}, we have
\begin{equation}\label{eq:ta}
 \frac{1}{2\pi} \TA(f^{\#}) 
     = -\chi(\overline{M}_{\gamma})+\sum_{j=1}^n(\mu_j^{\#}-d_j)
     = 2\gamma-2+\sum_{j=1}^n(\mu_j^{\#}-d_j)\;.
\end{equation}
Completeness of the metric $ds^2{}^{\#}$ at $p_j$ implies
 $\mu_j^{\#}-d_j\geq 1$. 
  However, the case $\mu_j^{\#}-d_j=1$ cannot occur (\cite[Lemma~3]{uy5}), so 
\begin{equation}\label{eq:complete}
   \mu_j^{\#}-d_j \geq 2 \; . 
\end{equation}
  
\subsection*{Effects of transforming the lift $\bmath{F}$.} 
Here we consider the change $\hat F=aFb^{-1}$ of the lift $F$,
where $a,b\in \SL(2,\C)$.  Then $\hat F$ is also a holomorphic null
immersion, and the hyperbolic Gauss map $\hat G$, the secondary Gauss
map $\hat g$ and the Hopf differential $\hat Q$ of $\hat F$ 
are given by (see \cite{uy3})
\begin{equation}\label{eq:changes}
  \hat G=a \star G, \quad \hat g=b\star g,\quad \hat Q=Q\;.
\end{equation}
In particular, the change $\hat F=aF$ moves the surface by a rigid
motion of $H^3$, and does not change $g$ and $Q$.  By choosing a 
suitable rigid motion $a\in \SL(2,\C)$ of 
the surface in $H^3$, we shall frequently use the following change 
of the hyperbolic Gauss map to  simplify its expression:
\begin{equation} \label{eq:simplifyG}
  \hat G =a\star G=\frac{a_{11}G+a_{12}}{a_{21}G+a_{22}} \; , 
  \qquad (a_{ij})_{i,j=1,2} \in \SL(2,\C) \; . 
\end{equation}
  
\subsection*{The Schwarzian derivative relation.} A direct computation implies 
that the secondary Gauss map $g$ depends on $G$ and $Q$ as
follows (\cite{uy1}):
\begin{equation}\label{eq:schwarz}
    S(g)-S(G)= 2Q\;,
\end{equation}
where
\[
   S(g)= \left[\left(\frac{g''}{g'}\right)' -
         \frac{1}{2}\left(\frac{g''}{g'}\right)^2\right]\,dz^2\qquad
   \left('=\frac{d}{dz}\right)
\]
is the {\em Schwarzian derivative\/} of $g$. 
Here, $z$ is a complex coordinate of $\overline{M}_\gamma$.

\subsection*{\boldmath$\SU(2)$-monodromy conditions.}  

Here we recall from \cite{ruy1} the construction of \cmcone{}
surfaces with given hyperbolic Gauss map $G$ and Hopf 
differential $Q$, which will play a crucial role in this paper.  
Let $\overline{M}_{\gamma}$ be a compact Riemann surface
and $M:=\overline{M}_{\gamma} \setminus\{p_1,\dots,p_n\}$.
Let $G$ and $Q$ be a meromorphic function and 
meromorphic $2$-differential on $\overline{M}_{\gamma}$.
The pair $(G,Q)$ must satisfy the following two compatibility
conditions:
\begin{align}
  \label{eq:cond1}
      &\text{For all $q \in M$, 
        $\ord_q Q$ is equal to the branching order of $G$, and}\\
  \label{eq:cond2} 
      &\text{for each end $p_j$, $\mu_j^{\#}-d_j\geq 2$.}
\end{align}
The first condition implies that the metric
${ds^2}^\#$ is (and hence ${ds^2}$ is also, by Lemma \ref{lem:complete}) 
non-degenerate at $q\in M$.
The second condition implies that the metric ${ds^2}^\#$ is complete (and 
hence ${ds^2}$ is also, again by Lemma \ref{lem:complete}) at 
$p_j\in \overline{M}_{\gamma}$ ($j=1,\dots,n$).

For such a pair $(G,Q)$, a solution $g$ of equation \eqref{eq:schwarz} has 
singularities at the branch points of $G$ (umbilic points or ends)
and the poles of $Q$ (ends).  However, regardless of whether $q \in M$ is 
a regular or umbilic point, $ds^2{}^{\#}$ and $Q^{\#}$ as in 
\eqref{eq:fund-forms} give a (non-degenerate) Riemann metric and holomorphic 
$2$-differential in a neighborhood $U_q \subset M$ of $q$.  
Then, by the fundamental theorem of surfaces, there exists a \cmcone{}
immersion $f^{\#}$ of $U_q$ into $H^3$ with induced
metric $ds^2{}^{\#}$ and Hopf differential $Q^{\#}$.
So the hyperbolic Gauss map $g$ of $f^{\#}$, which is a
solution of \eqref{eq:schwarz}, is a well-defined meromorphic function 
on $U_q$.
Since the solution of \eqref{eq:schwarz} is unique up to M\"obius
transformations $g\mapsto a\star g$ ($a\in\SL(2,\C)$), for 
any solution $g$ of \eqref{eq:schwarz} defined on the universal 
cover $\widetilde M$ of $M$, there exists a representation 
\[ 
    \rho_g\colon{}\pi_1(M)\to\PSL(2) \quad\text{ such that }\quad
    g\circ\tau^{-1} = 
    \rho_g(\tau)\star g 
\] 
for each covering transformation $\tau\in\pi_1(M)$. 

We now consider when the dual $f=(f^{\#})^{\#}$ (with data $(G,Q)$) 
of $f^{\#}$ is well-defined 
on $M$.  Choosing $F$ so that $F^{-1}$ is a lift of $f^{\#}$ (and then also 
$(F^{-1})^{-1}=F$ is a lift of $(f^{\#})^{\#}=f$), 
and noting that the representation $\rho_g\colon{}\pi_1(M)\to \PSL(2,\C)$ can 
be lifted into $\SL(2,\C)$ \cite{ruy1}, \eqref{eq:changes} implies 
\begin{equation}\label{eq:su2condition}
F^{-1} \circ \tau^{-1} = \rho_g(\tau) F^{-1} \end{equation} 
for each $\tau \in \pi_1(M)$.  Thus 
\begin{equation}\label{eq:su2condition2} 
f\circ \tau^{-1} = (F\circ\tau^{-1})(F\circ\tau^{-1})^{*} = 
F (\rho_g(\tau))^{-1} ((\rho_g(\tau))^{-1})^{*} F^{*} \; , \end{equation} 
and so $f$ is well-defined on $M$ if $\rho_g(\tau) \in \SU(2)$ for all 
$\tau \in \pi_1(M)$.  This is the crux of the following
Lemma~\ref{lem:exist}.  
Before stating this lemma, we need a definition:  

\begin{defn}
  A \cmcone{} immersion $f\colon{}M\to H^3$ is {\it  reducible\/} if 
  $\{\rho_g(\tau)\}_{\tau\in \pi_1(M)}$ are simultaneously 
  diagonalizable (i.e. if there exists a $P \in \PSL(2,\C)$ such that 
  $P\rho_g(\tau) P^{-1}$ is diagonal for all $\tau\in \pi_1(M)$).  
  If $f$ is not reducible, it is called {\it irreducible}.  
  When $f$ is reducible, it is either {\it $\Hyp^3$-reducible\/} or 
  {\it $\Hyp^1$-reducible\/} \cite{ruy1}, 
  and $f$ is called $\Hyp^3$-reducible if $\{\rho_g(\tau)\}_{\tau\in \pi_1(M)}$
  are all the identity, and is called $\Hyp^1$-reducible otherwise.  
\end{defn}

Clearly $f$ is $\Hyp^3$-reducible if and only if 
the lift $F$ itself is single-valued on $M$, by \eqref{eq:su2condition}.  
The name $\Hyp^1$-reducibility (resp.~$\Hyp^3$-reducibility) comes from
the fact that the surface has exactly a $1$ (resp.~$3$) dimensional 
deformation through surfaces preserving $G$ and $Q$ and the mean curvature, 
which is identified with the $1$ (resp.~$3$) dimensional hyperbolic space 
$\Hyp^1$ (resp.~$\Hyp^3$) \cite{ruy1}. 
On the other hand, if $f$ is irreducible, $f$ has no deformation preserving
mean curvature and $(G,Q)$ (see \cite{uy3,ruy1}).  

\begin{lemma}[\cite{uy3}]\label{lem:exist}
 Let $G$ and $Q$ be a meromorphic function and a meromorphic 
 $2$-differential on $\overline{M}_\gamma$ satisfying \eqref{eq:cond1} and 
 \eqref{eq:cond2}.  Assume $g$ is a solution of \eqref{eq:schwarz} such that 
 the image of $\rho_g$ lies in $\PSU(2)$.
 Then there exists a complete \cmcone{} immersion $f:M \to H^3$ with 
 hyperbolic Gauss map $G$, Hopf differential $Q$, and secondary Gauss map $g$. 
 
 If $f$ is irreducible, then $f$ is the unique surface with data 
 $(G,Q)$.  If $f$ is $\Hyp^1$-reducible $($resp.~$\Hyp^3$-reducible$)$, then 
 there exists exactly a $1$ $($resp.~$3)$ parameter family of \cmcone{} surfaces 
 with data $(G,Q)$.  
\end{lemma}

In the case that $M$ is of genus $\gamma=0$ with at most two ends, $f$ is 
reducible, as the fundamental group is commutative.  More generally, for the 
case $\gamma=0$ with $n$ ends, 
by Lemma~\ref{lem:exist} and the theory of linear ordinary differential
equations (see Appendix~\ref{app:log}), we have: 

\begin{proposition}\label{prop:exist}
 Let $\overline{M}_0=\C\cup\{\infty\}$ and  
 $M=\overline{M}_0\setminus\{p_1,\dots,p_n\}$ with $p_1,\dots,p_{n-1}$ 
 $\in \C$.
 Let $G$ and $Q$ be a meromorphic function and a meromorphic
 $2$-differential on $\C\cup\{\infty\}$ satisfying \eqref{eq:cond1} 
 and \eqref{eq:cond2}. 
 Consider the linear ordinary differential equation 
\begin{equation}\tag{E.0}
\label{eq:ode}
    \frac{d^2 u }{dz^2}+ r(z)u = 0 \; ,
\end{equation}
where $r(z)\,dz^2:=(S(G)/2)+Q$.  Suppose $n\ge 2$, and also 
$d_j=\ord_{p_j}Q\ge -2$ and the indicial equation of \eqref{eq:ode} 
at $z=p_j$ has the two roots $\lambda_1^{(j)},\lambda_2^{(j)}$ and 
log-term coefficient $c_j$, for $j=1,2,\dots,n-1$.  
\begin{enumerate}
 \item\label{item:prop:exist:1}
       Suppose that $\lambda_1^{(j)}-\lambda_2^{(j)}\in \Z^+$ and $c_j=0$ for 
       $j \leq n-1$.  Then there is exactly a $3$-parameter family of 
       complete conformal \cmcone{} immersions of $M$ into $H^3$ with
       hyperbolic Gauss map $G$ and  Hopf differential $Q$. 
       Moreover, such surfaces are $\Hyp^3$-reducible.
 \item\label{item:prop:exist:2} 
       Suppose that $\lambda_1^{(j)}-\lambda_2^{(j)}\in \Z^+$ and $c_j=0$ for 
       $j \leq n-2$, and that $\lambda_1^{(n-1)}-\lambda_2^{(n-1)}\in \R 
       \setminus \Z$.  
       Then there exists exactly a $1$-parameter family of complete 
       conformal \cmcone{} immersions of $M$ into $H^3$ with hyperbolic 
       Gauss map $G$ and  Hopf differential $Q$. 
       Moreover, such surfaces are $\Hyp^1$-reducible.
\end{enumerate}
Here we denoted by $\Z^+$  the set of positive integers.
\end{proposition}
The ordinary differential equation \eqref{eq:ode}
has also been applied in \cite{mu} for constructing certain 
classes of $\Hyp^3$-reducible \cmcone{} surfaces.  

\begin{proof}
The general theory of Schwarzian derivatives shows (\cite[Chapter 4]{Yoshida}) 
that for a linearly independent pair $u_1$, $u_2$ of solutions of
\eqref{eq:ode}, the function $g:=u_1/u_2$ satisfies \eqref{eq:schwarz}.  
Conversely, any function $g$ satisfying $S(g)=r(z)\,dz^2$ is obtained in 
this way.  

If $\lambda_1^{(j)}-\lambda_2^{(j)}=m\in \Z^+$ and $c_j=0$, 
then there is a fundamental system of solutions of \eqref{eq:ode} in a 
neighborhood of $p_j$ of the form 
\begin{equation}\label{eq:fund-int}
	  u_1  = (z-p_j)^{\lambda_1^{(j)}}\varphi_1(z)\;,
	   \qquad
	  u_2  = (z-p_j)^{\lambda_1^{(j)}-m}\varphi_2(z)\;,
\end{equation}
where $\varphi_1(z)$ and $\varphi_2(z)$ are holomorphic and nonzero at
$z=p_j$.
Then $g:=u_1/u_2$ satisfies 
\begin{equation}\label{totheremark}
        g\circ\tau_j^{-1} = \begin{bmatrix}
			     1 & 0 \\
                   0 & 1
                          \end{bmatrix}\star g \; ,
\end{equation}
where $\tau_j$ is the covering transformation which corresponds to a small
loop around $z=p_j$, implying $\rho_g(\tau_j)=\mbox{identity}$.  
So for case \ref{item:prop:exist:1}, we have
$\rho_g(\tau_j)=\mbox{identity}$ for all $j=1,\dots,n-1$, and therefore also for 
$j=n$, which implies that $g$ is a meromorphic function on $\C\cup\{\infty\}$.
By Lemma~\ref{lem:exist}, there exists a conformal \cmcone{} immersion
$f_a$ on $M$ with the secondary Gauss map
$a\star g$ for all $a\in \SL(2,\C)$. 
If $a\in\SU(2)$, then $f_a$ coincides with $f_{\id}$ by 
\eqref{eq:su2condition2}, so we have that 
the $3$-parameter family $(f_{[a]})_{[a] \in 
\SL(2,\subC)/\SU(2)}$ are complete 
conformal \cmcone{} immersions with hyperbolic Gauss map $G$ and  Hopf 
differential $Q$.

We remark here that if $\lambda_1^{(j)}-\lambda_2^{(j)}=m \in \Z^+$ 
and $c_j\neq 0$, then the monodromy matrix $\rho_g(\tau_j)$ defined by
$g\circ\tau_j^{-1}=\rho_g(\tau_j) \star g$ is not diagonalizable and 
is not even in $\SU(2)$.  So any \cmcone{} immersion on $\widetilde{M}$ 
(with $G$ and $Q$) cannot be well-defined on $M$ when some $c_j \neq 0$.  

Next we consider case \ref{item:prop:exist:2}, that is
$\lambda_1^{(n-1)}-\lambda_2^{(n-1)}\not\in \Z$.
There exists a fundamental system of
solutions of \eqref{eq:ode} of the form
\begin{equation}\label{eq:fund-no-int}
	  u_1  = (z-p_{n-1})^{\lambda_1^{(n-1)}}\varphi_1(z)\;,
	   \qquad
	  u_2  = (z-p_{n-1})^{\lambda_2^{(n-1)}}\varphi_2(z)\;,
\end{equation}
where $\varphi_1(z)$ and $\varphi_2(z)$
are holomorphic and nonzero at $z=p_{n-1}$.  
When $\tau_{n-1}$ is the covering transformation induced from a small loop 
about $z=p_{n-1}$, $g:=u_1/u_2$ satisfies
\begin{equation}
        g\circ\tau_{n-1}^{-1} = 
        \begin{bmatrix}
	  e^{\pi i (\lambda_1^{(n-1)}-\lambda_2^{(n-1)})} & 0 \\
          0 & e^{\pi i (\lambda_2^{(n-1)}-\lambda_1^{(n-1)})}
        \end{bmatrix}\star g \; .  
\end{equation}
In particular, $\rho_g(\tau_{n-1})\in\SU(2)$.  
On the other hand, in the proof of \ref{item:prop:exist:1}, we have seen
that $\rho_g(\tau_j)=\mbox{identity}$ for $j \in (1,\dots,n-2)$.
Hence $\rho_g(\tau_j)\in\SU(2)$ and are diagonal matrices 
for all $j \in (1,\dots,n)$, and we are in the $\Hyp^1$-reducible case.  
Note that this remains true when 
$g$ is replaced by 
\[
    s g(z)=a(s) \star g \; , 
      \quad \text{where} \quad 
     a(s) := 
        \begin{bmatrix}
	  \sqrt{s} & 0 \\
          0 & 1/\sqrt{s}
        \end{bmatrix} \; , 
     \quad \text{with} \quad s \in \R^+ \; , 
\]  
where $\R^+$ is the set of positive reals.
So we have a one-parameter 
family of complete conformal \cmcone{} immersions with hyperbolic Gauss map 
$G$ and Hopf differential $Q$ and secondary Gauss maps $s g$ for 
$s \in \R^+$.  ($s_1 g$ and $s_2 g$  for $s_1 \neq 
s_2$ will not produce equivalent surfaces, as $a(s_1)(a(s_2))^{-1} \not\in \SU(2)$.)  
Furthermore, Lemma~\ref{lem:exist} implies there is {\it only\/} 
a one-parameter family of \cmcone{} immersions with data $(G,Q)$.  
\end{proof}

By \eqref{eq:Fsharp}, we have 
\[
    (F^{-1})^{-1}d(F^{-1}) = 
          \begin{pmatrix} 
             g^\# & -g^\#{}^2 \\ 1 & -g^\#\hphantom{^2} 
          \end{pmatrix}
          \omega^\#\;,
\]
where
\[
      g^\#=G \; ,\qquad \omega^\#=-\frac{Q}{dG} \; .
\]
By Lemma~2.1 of \cite{uy1} (replacing $F$ with $F^{-1}$), we have that 
$X=F_{21}(z),F_{22}(z)$ satisfies the equation
\begin{equation}\tag*{(E.1)$^{\#}$}
\label{eq:E-1}
   X''-\bigl(\log(\hat{\omega}^\#)\bigr)'X'+\hat Q X = 0 \; , 
\end{equation}
and $Y=F_{11}(z),F_{12}(z)$ satisfies the equation
\begin{equation}\tag*{(E.2)$^{\#}$}
\label{eq:E-2}
   Y''-\bigl(\log(G^2 \hat{\omega}^\#)\bigr)'Y'+\hat Q Y = 0 \; , 
\end{equation}
where $Q(z)=\hat Q(z) dz^2$ and $\omega^\#= \hat{\omega}^\# (z) dz$. (We call 
them \ref{eq:E-1} and \ref{eq:E-2} because they are the dual versions of 
equations (E.1) and (E.2) in \cite{uy1}.)  
These two equations have been shown in \cite{Yu2} as a modification
of the corresponding equations in \cite{uy1}.
As we will see later, equations \ref{eq:E-1} and \ref{eq:E-2} are 
sometimes more convenient than equation \eqref{eq:ode} for solving 
monodromy problems.  In fact, we will have use for the following lemma: 
\begin{lemma}\label{lem:added} 
 Let $G$ and $Q$ be a meromorphic function and a holomorphic 
 $2$-dif\-ferential on $D^*=\{z\in \C\,;\, 0<|z|<1 \}$ such that 
 the metric $ds^2{}^\#$ defined by \eqref{eq:fund-forms} 
 is positive definite on $D^*$ and complete at $0$.  
 Assume $\ord_{z=0} Q \geq -2$ and $Q$ is not identically zero.  
 Then the following three conditions are all equivalent.
\begin{enumerate}
\item\label{item:added:1}
      The difference of the solutions of the indicial equation of \/
      {\rm \ref{eq:E-1}} at $z=0$ is a positive integer and the log-term 
      coefficient of {\rm \ref{eq:E-1}} vanishes.
\item\label{item:added:2} 
      The difference of the solutions of the indicial equation of 
      {\rm \ref{eq:E-2}} at $z=0$ is a positive integer and the log-term 
      coefficient of {\rm \ref{eq:E-2}} vanishes.
\item\label{item:added:3}
      The difference of the solutions of the indicial equation of 
      \eqref{eq:ode} at $z=0$ is a positive integer and the log-term 
      coefficient of \eqref{eq:ode} vanishes.  
\end{enumerate}
\end{lemma}

\begin{proof}
The hyperbolic Gauss map of the dual surface $f^\#=F^{-1}(F^{-1})^*$
is equal to the secondary Gauss map $g$ of $f=F F^*$.
Thus conditions \ref{item:added:1} and \ref{item:added:2} are equivalent to 
the condition that $g$ is single valued at $z=0$, by Lemma~2.2 
of \cite{uy1}. On the other hand, as seen in the proof of 
Proposition~\ref{prop:exist}, condition \ref{item:added:3} is also 
equivalent to the condition that $g$ is single valued at $z=0$.  
\end{proof}
Here is a natural place to include the next lemma, which we shall use in 
the sequel, \cite{ruy4}, to this paper.  
\begin{lemma}\label{lem:added2} 
 With the same assumptions as in Lemma~{\rm\ref{lem:added}}, 
 the following three conditions are all equivalent.  
\begin{enumerate}
\item\label{item:added2:1}
      The difference of the solutions of the indicial equation of 
      {\rm \ref{eq:E-1}} at $z=0$ is a real number.  
\item\label{item:added2:2}
      The difference of the solutions of the indicial equation of 
      {\rm \ref{eq:E-2}} at $z=0$ is a real number.  
\item\label{item:added2:3}
      The difference of the solutions of the indicial equation of 
      \eqref{eq:ode} at $z=0$ is a real number.  
\end{enumerate}
\end{lemma}
\begin{proof}
We write 
\[ 
    G(z) = z^\mu \hat{G}(z) \; , 
    \qquad
    \omega^\# (z) = z^\nu \hat{\omega}^\# (z) dz \; , 
\] 
where $\hat{G}$ and $\hat{\omega}^\#$ are 
nonzero and holomorphic at $z=0$, for some integers $\mu$ and $\nu$.  

If $\ord_{z=0} Q = -2$, so $\mu+\nu=-1$ and $Q = (\theta z^{-2}+\dots)dz^2$ 
for some $\theta \neq 0$, then the difference of the solutions of 
the indicial equations is $\sqrt{\mu^2-4 \theta}$ in all three cases, hence 
the three statements are clearly equivalent.  

If $\ord_{z=0} Q \geq -1$, then the indicial equation in the first case 
(resp.~second case, third case) is 
\[ 
     t(t-1)-\nu t = 0\;, \quad
     \left(
     \text{resp.}~
     t(t-1)-(2\mu + \nu ) t = 0 \; , 
     \quad
      t(t-1) + \frac{1-\mu^2}{4} =0 \;
     \right)
     \; . 
\]  
Hence the difference of the roots is $|\nu+1|$ 
(resp.~$|2\mu+\nu+1|$, $|\mu|$), and so all three statements hold.  
\end{proof}
\section{The classification of surfaces with $\TA(f^{\#})\le 4\pi$}
\label{sec:4-pi}

We begin our consideration of classification with this simple case:
\begin{theorem} 
  A complete \cmcone{} immersion $f$ with 
  $\TA(f^{\#})\le 4\pi$ is congruent to one of the following{\rm :}
\begin{enumerate}
 \item a horosphere, 
 \item an Enneper cousin dual,  
       $(g,Q)=(\tan\sqrt{\theta}z,\theta dz^2)$ 
       $(\theta\in \C\setminus \{0\})\;$, 
 \item a catenoid cousin,
\[
    (g,Q)=
       \left(a z^{\mu},\frac{1-\mu^2}{4z^2}\,dz^2\right)
      \qquad
      (a\in\R^+\;,\;\mu\in\R^+\setminus\{1\})\;,
\]
\item a warped catenoid cousin that has a degree $1$ hyperbolic
      Gauss map,
\[
    (g,Q)=
      \left(a z^l+b,\frac{1-l^2}{4z^2}dz^2\right)\qquad
         (a,b\in\C\setminus\{0\}\;,\; l \in \Z^+\setminus\{1\})\;.
\]
\end{enumerate}
\end{theorem}

\begin{proof}
  Since $\TA(f^{\#})  \in 4\pi\Z$, 
we need to consider only the cases $\TA(f^{\#})=0$ and $4\pi$. 
If $\TA(f^{\#})=0$, then the hyperbolic Gauss map is constant, so 
\eqref{eq:pseudo} implies $K^{\#}\equiv 0$.  Thus $f^{\#}$ is a 
totally umbilic \cmcone{} immersion, so 
both $f^{\#}$ and $f$ are horospheres.  
So we consider the remaining case $\TA(f^{\#})=4\pi$.  
Then $G$ is meromorphic of degree $1$ on $\overline M_{\gamma}$,
which implies $\gamma=0$.  Hence we may choose $\overline M_0=
\C\cup\{\infty\}$, and by \eqref{eq:simplifyG}, we may assume  $G=z$.  
Since $G$ has no branch points, \eqref{eq:QandG} implies 
there are no umbilic points, and \eqref{eq:mu-ineq} implies 
\begin{equation}\label{eq:mu-4pi}
   \mu_j^{\#}=0\
\end{equation}
at each end $p_j$.  By \eqref{eq:ta} and \eqref{eq:mu-4pi} and the fact that 
$\gamma=0$, we have
\begin{equation}\label{eq:ineq-4pi}
     2 =\frac{1}{2\pi} \TA(f^\#)=
      -2-\sum_{j=1}^{n} d_j \; .  
\end{equation}
By \eqref{eq:complete}, we have $2\ge -2+2n$, so $n=1$ or $2$.  

\subsection*{The case \boldmath$n=1$}  
In this case, \eqref{eq:ineq-4pi} implies $d_1=-4$.  We may put 
the end at $p_1=\infty$, and then $Q$ has a single pole of order $4$ at
$\infty$ and no zeroes.  
Thus $Q=\theta\,dz^2$ for some $\theta\in\C\setminus\{0\}$.

A \cmcone{} surface in $H^3$ with secondary Gauss map $g=z$
and Hopf differential $Q=\theta\,dz ^2$ is called an Enneper cousin
\cite{Bryant}.
So a surface with data $(G,Q)=(z,\theta \, dz ^2)$ is the dual
of an Enneper cousin \cite[Example~5.4]{ruy1}. 
(Recall that dualizing switches the two Gauss maps, and changes the Hopf 
differential only by a sign.)  
\subsection*{The case \boldmath$n=2$}
In this case, \eqref{eq:ineq-4pi} becomes $4 = -d_1-d_2$.
Then $d_j=-2$ ($j=1, 2$), by \eqref{eq:complete}.
Hence the immersion $f$ is a \cmcone{} surface of genus $0$ whose two ends 
must both be regular \cite{uy1}, and this type of surface is 
classified in \cite{uy1}.  In particular, $f$ is in the case $m=1$ of 
Theorem~6.2 in \cite{uy1}.  So the surface is either a {\em catenoid 
cousin\/} \cite[Example~2]{Bryant} or 
a {\em warped catenoid cousin\/} with embedded ends 
(the case $m=1$ in Theorem~6.2 in \cite{uy1}).  
\end{proof}

The warped catenoid cousins are described in detail in \cite{ruy4}.  
\section{Surfaces with $\TA(f^{\#})=8\pi$}
\label{sec:8-pi}

We now assume $f$ has $\TA(f^{\#})=8\pi$.  
Then, by \eqref{eq:ta} and \eqref{eq:complete}, 
\begin{equation}\label{eq:ineq-8pi}
    6 = 2 \gamma + \sum_{j=1}^n (\mu_j^{\#}-d_j)
       \geq 2(\gamma+n)
\end{equation}
holds.  Thus the possible cases are
\[
     (\gamma,n)  = (0,1)\;, \quad (0,2)\;, \quad (0,3)\;, \quad
                   (1,1)\;, \quad (1,2)\;, \quad \text{and}\quad
                   (2,1)\;.
\]
Since $\TA(f^{\#})=8\pi$, $G$ is meromorphic 
on $\overline M_{\gamma}$ of degree $2$.  Hence \eqref{eq:mu-ineq} implies 
\begin{align}\label{eq:degree}
    &\mu_j^{\#}\leq 1 \quad &(&j=1,2,\dots, n)\;,\\
\intertext{and at each umbilic point $q_l$,} 
             \label{eq:ord-umb}
    &\xi_l=1  \quad         &(&l=1,2,\dots,k)\;.
\end{align}

\subsection*{The case \boldmath$(\gamma,n)=(2,1)$}  
Since equality holds in \eqref{eq:osserman}, the single 
end $p_1$ is embedded.  
By \eqref{eq:ineq-8pi}, $\mu_1^{\#}-d_1=2$.  Thus the possible cases are
\[
   (\mu_1^{\#},d_1) = (0,-2) \quad \text{or}\quad (1,-1) \; , 
\]
by \eqref{eq:degree}.  
If $(\mu_1^{\#},d_1)=(0,-2)$, the end $p_1$ is of type I in the sense of
\cite{ruy2}, so the flux about this end does not vanish
\cite[Proposition~2]{ruy2}. 
If $(\mu_1^{\#},d_1)=(1,-1)$, then, since the end is embedded, Corollary~5 
in \cite{ruy2} implies that the flux about the end again does not
vanish.  But non-vanishing flux at a single end contradicts the balancing
formula \cite[Theorem~1]{ruy2}, so the case $(\gamma,n)=(2,1)$ does not
occur.  

\subsection*{The case \boldmath$(\gamma,n)=(1,2)$}  
In this case, \eqref{eq:ineq-8pi} implies $4 = (\mu_1^{\#}-d_1)+ (\mu_2^{\#}-
d_2)$.  By \eqref{eq:complete}, we have $\mu_j^{\#}-d_j = 2$ for 
$j=1,2$.  Hence \eqref{eq:degree} implies 
\[
   (\mu_j^{\#},d_j)=(0,-2)\quad\text{or}\quad (1,-1)\qquad
    (j=1,2)\;.
\]

Assume $d_1=-2$ and $d_2=-1$.
Then, by the transformation \eqref{eq:simplifyG} if necessary, 
we may assume the hyperbolic Gauss map has a zero or pole at each end.  
In this case, the end $p_1$ is regular of type I, and 
$p_2$ is regular of type II in the sense of \cite{ruy2}, 
contradicting Theorem~7 in \cite{ruy2}.  
Hence this case is impossible, leaving the two remaining possibilities:
\begin{align}
   (\mu_1^{\#},d_1)=(\mu_2^{\#},d_2)&=(0,-2)\;, \label{eq:I(-2,-2)}\\
   (\mu_1^{\#},d_1)=(\mu_2^{\#},d_2)&=(1,-1)\;. \label{eq:I(-1,-1)}
\end{align}
\begin{figure}
\begin{center}
\vspace*{0.1in}
\begin{tabular}{ccc}
  \includegraphics[width=1.5in]{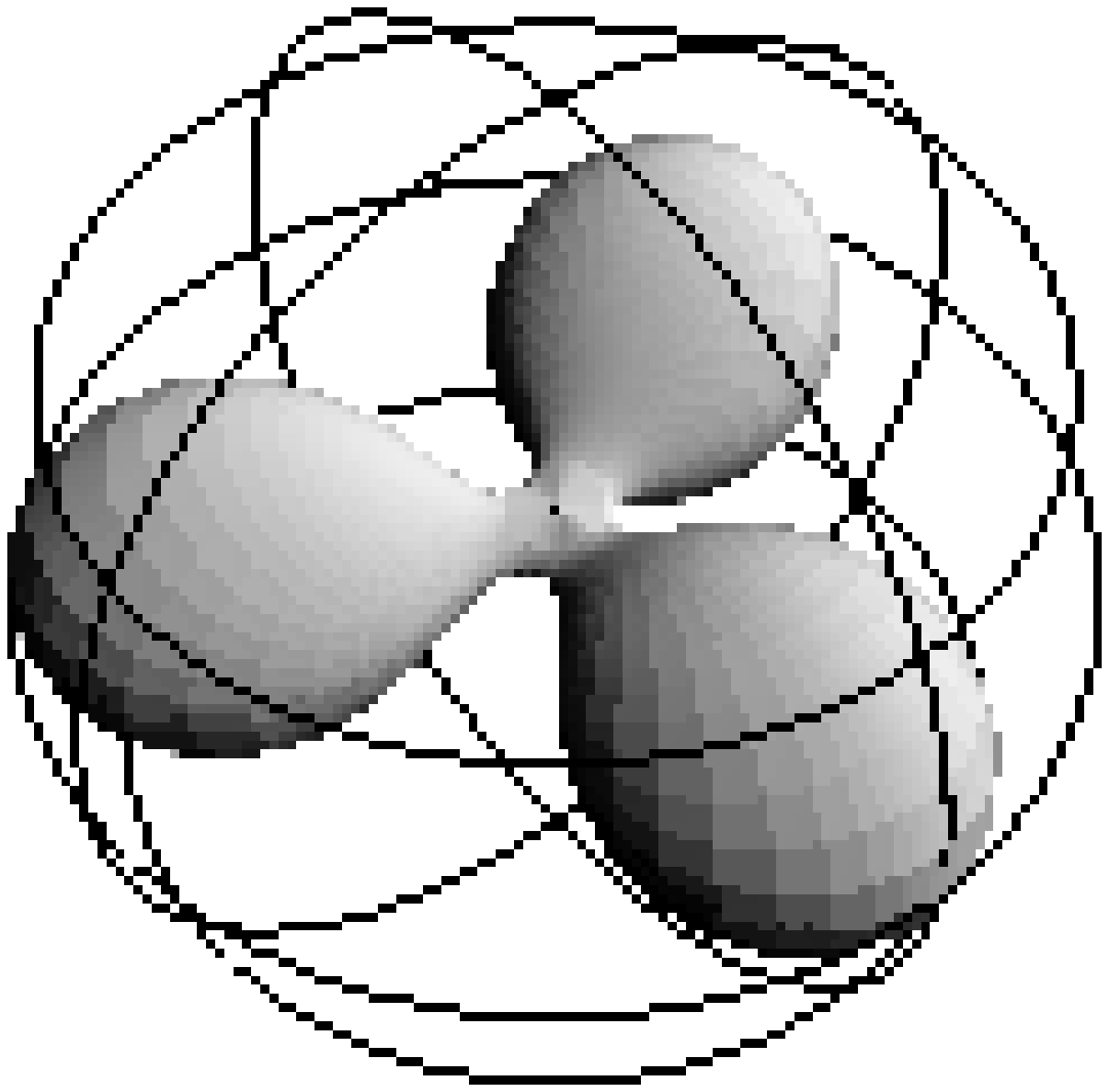} &
  \includegraphics[width=1.5in]{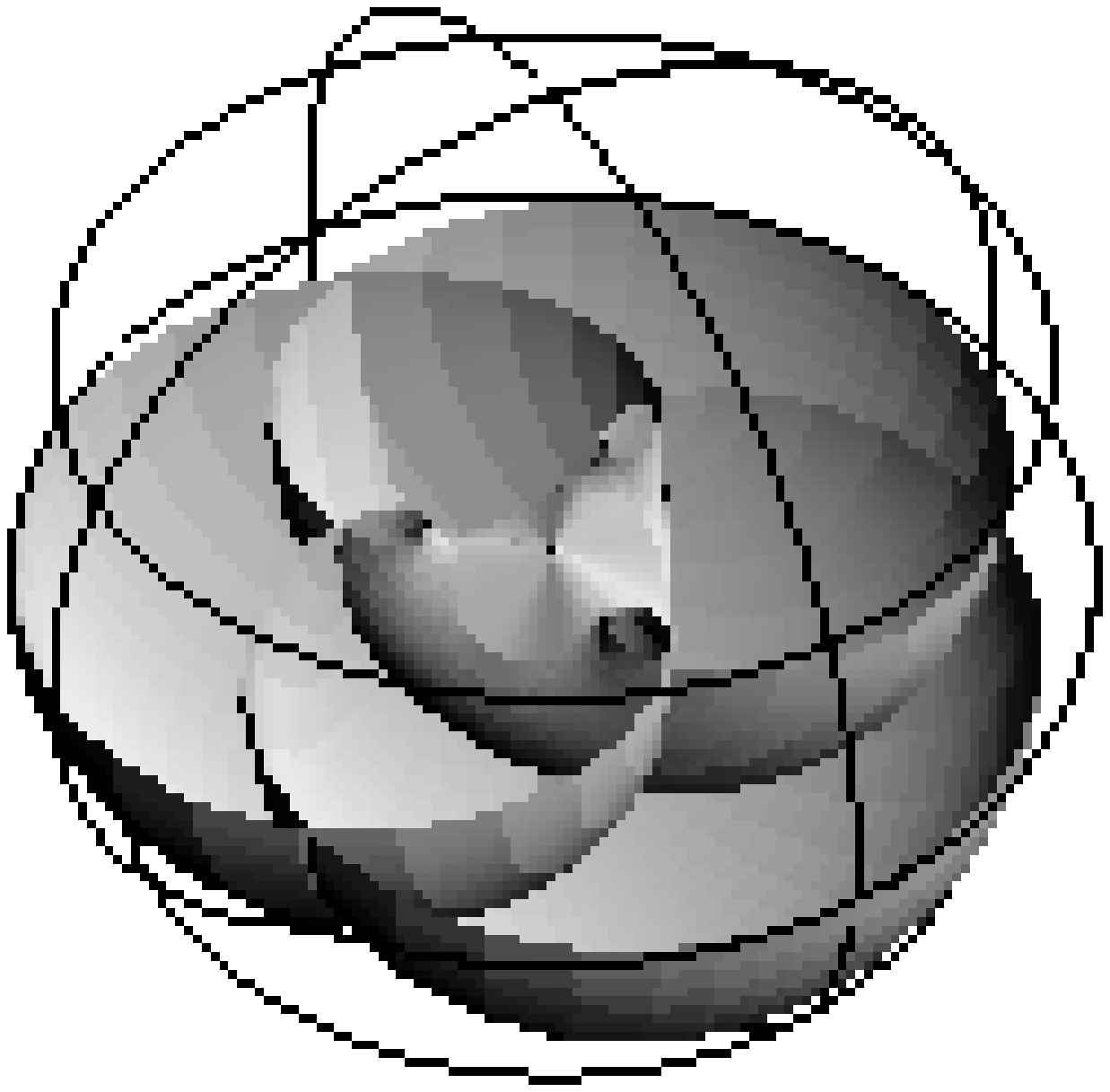} &
  \includegraphics[width=1.5in]{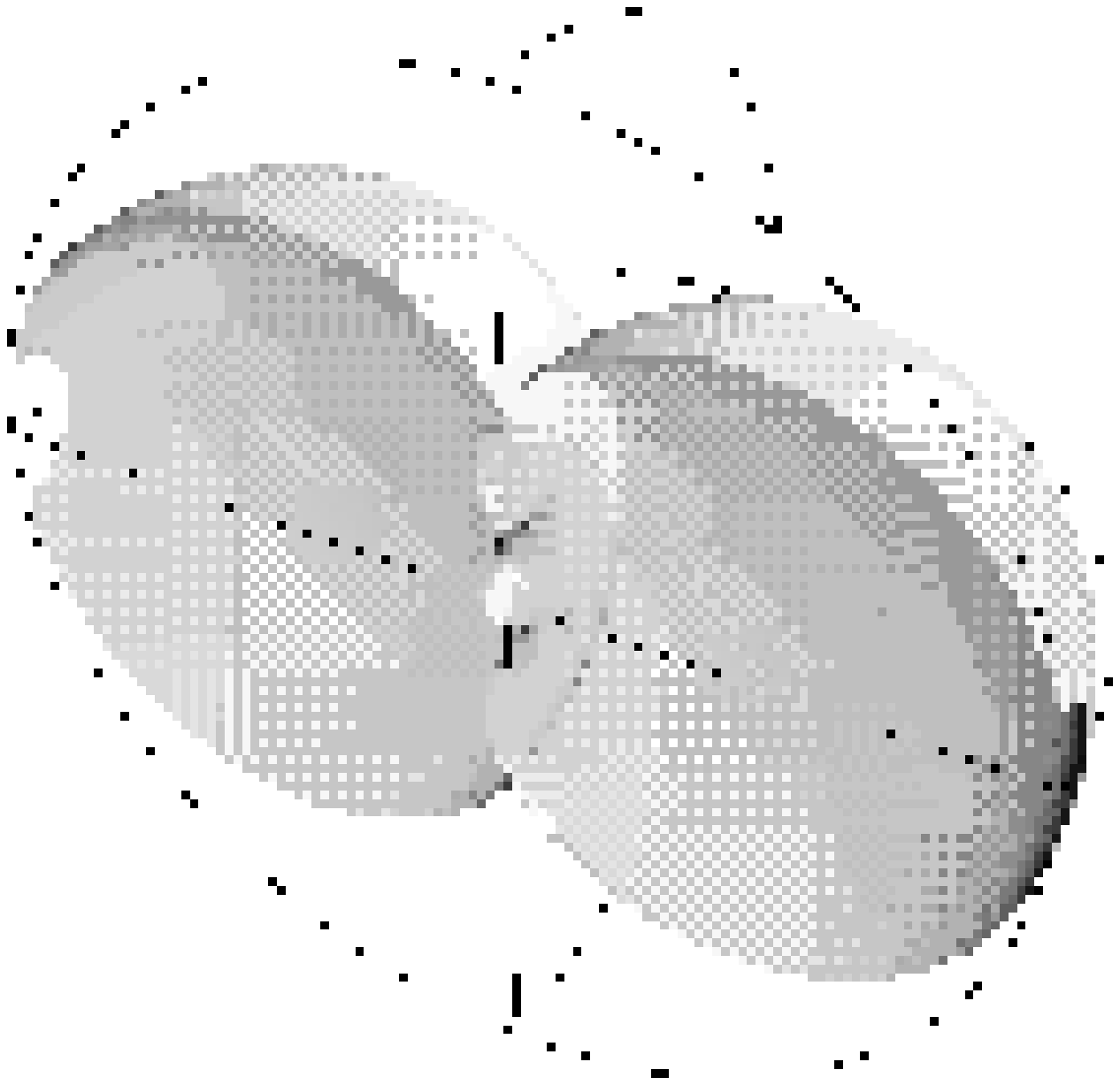} 
\end{tabular}
\vspace*{0.1in}
\end{center}
\caption{Two \cmcone{} trinoids in $H^3$, which are surfaces of type 
$\gO(-2,-2,-2)$, and a genus $1$ catenoid cousin, which is a surface 
of type $\gI(-2,-2)$, shown in the Poincar\'e model of $H^3$.  
Only one of two congruent pieces of the right-most two surfaces is 
shown, and the other half of each surface is the reflection (in the plane 
containing the boundary curves seen here) of the piece shown.}
\end{figure}
For the case \eqref{eq:I(-2,-2)}, the first author  and Sato \cite{rs}
constructed a one-parameter family of ``genus one catenoid cousins''.  
Note that such surfaces cannot exist as minimal 
surfaces in $\R^3$, by Schoen's result \cite{s}.  

\subsubsection*{Surfaces of type \boldmath{$\gI(-1,-1)$}}
For the case \eqref{eq:I(-1,-1)},
we can determine the candidates of $(G,Q)$ explicitly as follows
(however, the period problem is unsolved and no example is known):
\begin{proposition}\label{prop:data-I}
Let $\overline{M}_1=\C/\Gamma$, where $\Gamma$ is a lattice on $\C$, 
and assume there exists a \cmcone{} immersion
 $f\colon{}\overline{M}_1\setminus\{p_1,p_2\}\to H^3$ 
with\/ $\TA(f^{\#})=8\pi$ of type\/ $\gI(-1,-1)$.
Then there exists a generating pair $\{v_1,v_2\}\subset\C$ of\/ $\Gamma$
such that the hyperbolic Gauss map $G$ and Hopf
differential $Q$ are given by 
\begin{equation}\label{eq:data-I(-1,-1)}
   G=\wp(z)\;,\qquad 
   Q(z)=\theta \,
   \frac{\sigma(z-v_1/2)\sigma(z-v_2/2)}{\sigma(z)\sigma(z-(v_1+v_2)/2)}
   \,dz^2
   \qquad (\theta \in \C\setminus\{0\})\;,
\end{equation}
where $\wp(z)$ is the Weierstrass $\wp$-function 
and $\sigma$ is the entire function defined by
\[
   \sigma(z):=z\prod_{v\in \Gamma\setminus\{0\}}
   \left\{\left(1-\frac{z}{v}\right)
   e^{\frac{z}v+\frac{z^2}{2v^2}}\right\}\;.
\]
\end{proposition}
\begin{proof}
In this case, the hyperbolic Gauss map $G$ is of degree $2$.
Without loss of generality, we may assume that $z=0$ is an end of
the surface. Moreover, by \eqref{eq:simplifyG} we may assume that 
$z=0$ is a pole of $G$.  
As $z=0$ is a branch point of $G$ (since $\mu_j^\#=1$), $G$ has a pole of 
order $2$ at $z=0$.  Up to a constant multiple, the function $\wp(z)$ is 
uniquely characterized 
as a degree $2$  meromorphic function on $\C/\Gamma$ with a pole of
order $2$ at the origin \cite{hc}. Thus we have $G(z)=c\,\wp(z)$, and 
we can normalize $c=1$, by \eqref{eq:simplifyG}.

Suppose $\{v_1,v_2\}$ generates $\Gamma$. 
Then the branch points of $\wp$ are $0$, $v_1/2$, $v_2/2$ and $(v_1+v_2)/2$
modulo $\Gamma$, which are the ends and umbilic points.  
We assume $0$ and $(v_1+v_2)/2$ are the ends.  
(If $v_1/2$ is an end, for example, we may change the generator $\Gamma$ to 
$\{ \tilde v_1 = v_1-v_2,\tilde v_2=v_2 \} $.)  
Thus the umbilic points are $v_1/2$ and $v_2/2$.

Next we find the Hopf differential $Q(z)=q(z)\,dz^2$, using the following 
fact:  
\begin{fact}[\cite{hc}]
  Let $a_1,\dots,a_n$ and $b_1,\dots,b_n$ be points in $\C$ such
  that $a_j\ne b_k \pmod \Gamma$, $j,k \in \{1,\dots,n\}$, and 
  $\sum_{j=1}^n a_j = \sum_{k=1}^n b_k \pmod \Gamma$.  Then
\[
  f(z):=\theta\frac{\sigma(z-a_1)\cdots\sigma(z-a_n)}%
        {\sigma(z-b_1)\cdots\sigma(z-b_n)} \qquad 
       (\theta\in\C\setminus\{0\})
\]
  is a meromorphic function on $\C/\Gamma$ such that
  $\{a_1,\dots,a_n\}$ $($resp.~ $\{b_1,\dots,b_n\})$ are the
  set of zeroes $($resp.~poles$)$, i.e. the divisor of $f$ is 
  $a_1+\dots+a_n-b_1-\dots-b_n$.  
  Conversely, any elliptic function on $\C/\Gamma$ with the same 
  divisor is of this form.  
\end{fact}
The meromorphic function $q(z)$ should have poles of order $1$ at $z=0$,
$(v_1+v_2)/2$ (ends) and zeroes of order $1$ at $z=v_1/2$, $v_2/2$
(umbilic points).
Thus $Q(z)$ can be written as in \eqref{eq:data-I(-1,-1)}.
\end{proof}

\subsection*{The case \boldmath$(\gamma,n)=(1,1)$} 
By \eqref{eq:ineq-8pi} and \eqref{eq:degree}, we have two possible 
cases:
\[
    (\mu_1^{\#},d_1) = (0,-4) \quad \text{or} \quad (1,-3)\;.
\]
The second of these cases (the $\gI(-3)$ case) is still unknown, but 
for the first case $\gI(-4)$, the following proposition provides examples, 
proven (in Section~\ref{sec:deform}) 
by deforming from a complete minimal surface
in $\R^3$ of genus $1$ with one end satisfying $d_1=-4$.
\begin{proposition}\label{thm:gi4}
 By deforming the Chen-Gackstatter surface in $\R^3$ \cite{cg}, 
 one obtains a one-parameter family of \cmcone{} surfaces of
 type $\gI(-4)$ with dual total absolute curvature $8\pi$.  
\end{proposition}

\subsection*{The case \boldmath{$(\gamma,n)=(0,3)$}}  
Here, \eqref{eq:ineq-8pi} and \eqref{eq:complete} imply $\mu_j^{\#}-d_j=2$ 
for $j=1,2,3$.
Moreover, \eqref{eq:r-r} implies $d_1+d_2+d_3\leq -4$.  So 
\eqref{eq:degree} implies that the possibilities are:  
\begin{center}
\begin{tabular}{l@{ : }l}
 Type $\gO(-2,-2,-2)$ &
 $(d_1,d_2,d_3)=(-2,-2,-2)$ and $(\mu_1^\#,\mu_2^\#,\mu_3^\#)=(0,0,0)\;$,\\
 Type $\gO(-1,-2,-2)$ & 
 $(d_1,d_2,d_3)=(-1,-2,-2)$ and $(\mu_1^\#,\mu_2^\#,\mu_3^\#)=(1,0,0)\;$,\\
 Type $\gO(-1,-1,-2)$ &
 $(d_1,d_2,d_3)=(-1,-1,-2)$ and $(\mu_1^\#,\mu_2^\#,\mu_3^\#)=(1,1,0)\;$.
\end{tabular}
\end{center}
In each case, equality holds in \eqref{eq:osserman}, so all ends are 
embedded.  
Since the genus of the surface is $0$, we can set
$\overline{M}_0=\C\cup\{\infty\}$.

\subsubsection*{Surfaces of type \boldmath{$\gO(-2,-2,-2)$}}  
Such surfaces have three embedded ends with $d_j=-2$ ($j=1,2,3$), and the 
irreducible ones are classified in \cite[Theorem~2.6]{uy6}.  
So here we consider the reducible case. 

We may set $p_1=0$, $p_2=1$ and $p_3=\infty$.
By  \eqref{eq:r-r} and \eqref{eq:ord-umb}, there are two 
distinct umbilic points $q_1$ and $q_2$ of order $1$.
Then the Hopf differential $Q$ must have simple zeroes
at $q_1$ and $q_2$ and poles of order $2$ at $0$, $1$ and $\infty$.
Since all three $\mu_j^{\#}=0$, $q_1$ and $q_2$ are the only 
branch points of $G$.  
Also, $G(q_1)$, $G(q_2)$, and $G(\infty)$ are all distinct, because $q_1$ and 
$q_2$ are double points of $G$ and $\deg G=2$.  
Then, by \eqref{eq:simplifyG}, we can set $G(q_1)=0$, $G(q_2)=\infty$, 
and $G(\infty)=1$.  Thus $G$ and $Q$ are written as
\begin{equation}\label{eq:data-O(-2,-2,-2)}
  G = \left(\frac{z-q_1}{z-q_2}\right)^2,\qquad
  Q = \theta\,\frac{(z-q_1)(z-q_2)}{z^2(z-1)^2}\,dz^2\qquad
     (\theta\in\C\setminus\{0\})\;.
\end{equation}

\begin{example}[$\Hyp^1$-reducible examples of type $\gO(-2,-2,-2)$]
\label{ex:o222red1}
For $s \in \R$ such that
\begin{equation}\label{eq:o222dif2}
       -4\frac{1 + 4 s + s^2}{1 + 10 s + s^2}\in\R\setminus\Z\;,
\end{equation}
let 
\begin{equation}\label{eq:o222redq}
       q_1 = \frac{1+10s+s^2}{4s(1-s)}\;, \qquad
       q_2 = \frac{1+10s+s^2}{4(s-1)} \; , \quad
   \text{and}\quad
       \theta = -\frac{3}{4q_1q_2}\;.
\end{equation}
Consider \eqref{eq:ode} for $r(z)\,dz^2=(S(G)/2)+Q$, with $G$ and $Q$ 
determined by \eqref{eq:data-O(-2,-2,-2)} and \eqref{eq:o222redq}.
Then the roots of the indicial equation of
\eqref{eq:ode} at $z=0$ are $-1/2$ and $3/2$, so their difference is $2 
\in \Z$, and one can check by
\eqref{m=2} that the log-term coefficient vanishes.
Moreover, the difference of the roots of the indicial equation 
at $z=1$ equals the value in \eqref{eq:o222dif2}.
Hence, by \ref{item:prop:exist:2} of Proposition~\ref{prop:exist}, 
there exists an $\Hyp^1$-reducible \cmcone{} immersion 
$f\colon{}\C \setminus \{0,1\}\to H^3$ with $G$ and $Q$ as in
\eqref{eq:data-O(-2,-2,-2)} and \eqref{eq:o222redq}.
Since each surface is $\Hyp^1$-reducible (this follows from the fact that 
the difference of the roots of the indicial equation is an integer at 
$z=0$ and not an integer at $z=1$), there exists a one-parameter 
family of \cmcone{} surfaces for each $s$, with this $G$ and $Q$.  
Thus, we have found a $2$-parameter family of $\Hyp^1$-reducible \cmcone{}
surfaces of type $\gO(-2,-2,-2)$.
\end{example}
\begin{example}[$\Hyp^3$-reducible examples of type $\gO(-2,-2,-2)$]
       \label{ex:o222red3}
  For $m\geq 2$, $m \in \Z$, let 
\[
     q_1 = \frac{1}{2}\left(1+\frac{1}{\sqrt{m}}\right)\;,\qquad
     q_2 = \frac{1}{2}\left(1-\frac{1}{\sqrt{m}}\right)\;,\qquad
   \text{and}\quad
       \theta = -m(m+1) \; .  
\]
  Then a meromorphic function $g$ on $\C\cup\{\infty\}$ such that
\[
    dg = z^{m-1}(z-1)^{m-1}(z-q_1)(z-q_2) \,dz
\]
 satisfies equation \eqref{eq:schwarz} for $G$ and $Q$ as in 
\eqref{eq:data-O(-2,-2,-2)}.  Since $g$ is meromorphic, 
$\rho_g(\tau)$ is the 
identity for all $\tau \in \pi_1(\C \setminus \{0,1\})$, so 
Lemma~\ref{lem:exist} implies there exists an $\Hyp^3$-reducible 
\cmcone{} immersion 
 $f\colon{}\C\setminus\{0,1\}\to H^3$ whose hyperbolic Gauss map, 
 Hopf differential, and secondary Gauss map are $G$, $Q$, and $g$,
 respectively. 
\end{example}

\subsubsection*{Surfaces of type \boldmath{$\gO(-1,-2,-2)$}}
In this case, we will see that there is a $2$-parameter family of 
$\Hyp^1$-reducible surfaces, and countably many $\Hyp^3$-reducible
families.
By \eqref{eq:r-r}, there exists one umbilic point of order $1$.
Without loss of generality, we can set the ends to be 
$(p_1,p_2,p_3)=(0,1,p)$ ($p\in\C\setminus\{0,1\}$)
and the umbilic point to be $q_1=\infty$. 
Then the Hopf differential $Q$ has a pole of order $2$ (resp.~order $1$) 
at $z=1,p$ (resp.~$z=0$)
and has no zeroes on $\C$, so it has the form 
\[
  Q=\frac{\theta\, dz^2}{z(z-1)^2(z-p)^2}\qquad
   (\theta\in\C\setminus\{0\})\;.
\]
By \eqref{eq:QandG} and the fact $\mu_1^\#=1$, $G$ has branch points of
order $1$ at $z=0$ and $\infty$. 
Then, by \eqref{eq:simplifyG}, we may assume $G = z^2$, because $\deg G=2$.
Consider the ordinary differential equation \eqref{eq:ode} with 
$r(z)\,dz^2 = (S(G)/2)+Q$.  At the singularity $z=0$, $r(z)$ expands as 
\[
   r(z)= -\frac{3}{4}\frac{1}{z^2}+\frac{\theta}{p^2}\frac{1}{z}
            + \frac{2\theta(p+1)}{p^3}+O(z)\;.
\]
Thus the difference of the roots of the indicial equation of 
\eqref{eq:ode} at $z=0$ is $2$. 
Then, by \eqref{m=2}, the log-term coefficient of \eqref{eq:ode} at
$z=0$ vanishes if and only if $\theta=-2p(p+1)$.
Hence, if such a surface exists, $G$ and $Q$ are 
\begin{equation}\label{eq:data-O(-1,-2,-2)}
    G = z^2\; ,\qquad 
    Q = \frac{-2p(p+1)}{z(z-1)^2(z-p)^2}\,dz^2\qquad
    (p\in\C\setminus\{0,1\})\;.
\end{equation}
For $G$ and $Q$ as in \eqref{eq:data-O(-1,-2,-2)}, $r(z)$ expands at the 
singularity $z=1$ as 
\[
  r(z)= \frac{-2p(p+1)}{(1-p)^2}\frac{1}{(z-1)^2}+O\bigl((z-1)^{-1}\bigr)\;.
\]
Then the roots of the indicial equation of 
\eqref{eq:ode} at $z=1$ are 
\[
    \lambda_1 = 2 + \frac{2}{p-1}\;, \qquad
    \lambda_2 = -1 - \frac{2}{p-1}\;.
\]
So $\lambda_1-\lambda_2 \in \Z$ exactly when $4/(p-1) \in \Z$. 
Then, by Proposition~\ref{prop:exist}, we have
\begin{theorem}\label{thm:go122}
 Let $p \in \R$ such that $p \neq 1$ and $4/(p-1) \not\in \Z$.  
 Then there exists a conformal $\Hyp^1$-reducible \cmcone{} immersion
 $f\colon{}M=\C\cup\{\infty\}\setminus\{0,1,p\}$ with
 $\TA(f^{\#})=8\pi$ and hyperbolic Gauss map and Hopf differential as 
 in \eqref{eq:data-O(-1,-2,-2)}. 
 Moreover, all $\Hyp^1$-reducible surfaces with $\TA(f^{\#})=8\pi$ 
 of type $\gO(-1,-2,-2)$ are given in this manner.
\end{theorem}

The above discussion yields that all \cmcone{} surfaces of 
type $\gO(-1,-2,-2)$ are reducible.  So it only remains to 
classify the $\Hyp^3$-reducible case: 
\begin{theorem}
\label{thm:go122-h3}
 Let $r \geq 3$ be an integer and $p=(r+2)/(r-2)$.  
 Then there exists a conformal $\Hyp^3$-reducible \cmcone{} immersion
 $f\colon{}M=\C\cup\{\infty\}\setminus\{0,1,p\}$ with
 $\TA(f^{\#})=8\pi$ whose hyperbolic  Gauss map and Hopf differential are
 as in \eqref{eq:data-O(-1,-2,-2)}.
 Moreover, all $\Hyp^3$-reducible surfaces with $\TA(f^{\#})=8\pi$ of type
 $\gO(-1,-2,-2)$  are given in this manner.
\end{theorem}
\begin{proof}
 For given $r \geq 3$, there is a meromorphic function $g$ on 
 $\C\cup\{\infty\}$ so that 
\begin{equation}\label{eq:o122dg}
   dg = \frac{z(z-p)^{r-2}} {(z-1)^{r+2}}\,dz\;,
\end{equation}
 since the right-hand side of \eqref{eq:o122dg} has no residue.  
 One can check that $S(g)-S(G)=2Q$ when $p=(r+2)/(r-2)$.  
 Hence, by Lemma~\ref{lem:exist}, there exists an $\Hyp^3$-reducible 
 \cmcone{} immersion 
 $f\colon{}\C\cup\{\infty\}\setminus\{0,1,p\}\to H^3$ with 
 $G$ and $Q$ as in 
 \eqref{eq:data-O(-1,-2,-2)} and secondary Gauss map $g$ satisfying 
 \eqref{eq:o122dg}.

 Conversely, let $f\colon{}\C\cup\{\infty\}\setminus\{0,1,p\}\to H^3$ 
 be an $\Hyp^3$-reducible \cmcone{} immersion of type $\gO(-1,-2,-2)$ with 
 $\TA(f^{\#})=8\pi$. Then $G$ and $Q$ are as in \eqref{eq:data-O(-1,-2,-2)}. 
 Let $m_2$ (resp.~$m_3$) be the difference of the roots of the indicial
 equation of \eqref{eq:ode} at $z=1$ (resp.~$z=p$) for such $G$ and $Q$.
 Then we have $m_2 = |3 + (4/(p-1))|$ and $m_3 = |1 + (4/(p-1))|$.  
 Since $f$ is $\Hyp^3$-reducible, $m_2$ and $m_3$ are positive integers (so 
 also $4/(p-1) \in \Z$).  
 We may assume $m_2 \geq m_3$. (If not, we can exchange the two 
 ends $p$ and $1$, by changing $p$ and $z$ to $1/p$ and $z/p$.  Using 
 \eqref{eq:simplifyG}, we see that \eqref{eq:data-O(-1,-2,-2)} is 
 unchanged.)  
 
 Suppose that $m_2=m_3=1$, then $g$ is not branched at both 
 $1$ and $p$.  Noting that the branching orders of $g$ and $G$ are equal at 
 any finite point of the surface (this follows from equation 
 \eqref{eq:schwarz}), 
 we see that $g$ has branch points of order $1$ at $0$ and $\infty$ and no 
 other branch points.  So $g$ has degree $2$ and $g= a \star z^2$ for some 
 $a \in \SL(2,\C)$ and so $Q=(1/2) (S(g)-S(G)) = 0$, which is impossible.  
 
 Thus $m_2 \geq 2$, and 
 it follows that $4/(p-1)$ is a positive integer.  By setting $r = 
 2+(4/(p-1))\geq 3$, we have 
\[
     m_2 = 3 + \frac{4}{p-1} = r + 1 \;,\quad
     m_3 = 1 + \frac{4}{p-1} = r - 1 \;,\quad
   \text{and}\qquad
     p = \frac{r+2}{r-2}\;.
\]
 Thus $G$ and $Q$ are as in \eqref{eq:data-O(-1,-2,-2)} with
 $p=(r+2)/(r-2)$.  
\end{proof}

\subsubsection*{Surfaces of type \boldmath{$\gO(-1,-1,-2)$}}  
In this case, by \eqref{eq:r-r}, the surface has no umbilic points.  
We set the ends $(p_1,p_2,p_3)=(0,1,\infty)$.
The Hopf differential is then 
\begin{equation}\label{eq:go112-Q}
      Q = \frac{\theta\,dz^2}{z(z-1)}\;,\qquad 
    (\theta\in\C\setminus\{0\})\;.
\end{equation}
The hyperbolic Gauss map $G$ is a meromorphic function on 
$\C\cup\{\infty\}$ of degree $2$ with branch points of order $1$ at 
$z=0$ and $z=1$.  Hence we may set
\begin{equation}\label{eq:go112-G}
    G = \left(\frac{z-1}{z}\right)^2.
\end{equation}

\begin{theorem}\label{thm:go112}
 Any complete \cmcone{} immersion that is of type 
 $\gO(-1,-1,-2)$ with $\TA(f^{\#})=8\pi$ 
 is congruent to an $\Hyp^3$-reducible \cmcone{} immersion 
 $f\colon{}M=\C\setminus\{0,1\} \longrightarrow  H^3$ 
 with hyperbolic Gauss map and Hopf differential 
\[
   G = \left(\frac{z-1}{z}\right)^2\;,\qquad
   Q = \frac{-2\,dz^2}{z(z-1)} \; .
\]
\end{theorem}
\begin{proof}
 Consider equation \eqref{eq:ode}  for 
 $G$ and $Q$ in \eqref{eq:go112-G} and \eqref{eq:go112-Q} respectively.
 Then the roots of the indicial equations of 
 \eqref{eq:ode} are $-1/2$ and $3/2$ at both $z=0$ and $z=1$.  
 By \eqref{m=2}, the log-term coefficients at $z=0$ and 
 at $z=1$ both vanish if and only if $\theta=-2$.
 By Proposition \ref{prop:exist}, the corresponding $3$-parameter family of
 \cmcone{} immersions consists of immersions that are all well-defined
 on $M=\C\setminus\{0,1\}$  and are $\Hyp^3$-reducible.  
\end{proof}

\subsection*{The case \boldmath $(\gamma,n)=(0,2)$}  
In this case, \eqref{eq:ineq-8pi} and \eqref{eq:complete}
imply that 
\[
   (\mu_1^{\#}-d_1,\mu_2^{\#}-d_2) = (2,4)\quad\text{or}\quad
   (\mu_1^{\#}-d_1,\mu_2^{\#}-d_2) = (3,3)\;.
\]
Then, by \eqref{eq:degree}, all possibilities are:  
\begin{center}
\begin{tabular}{l@{ : }l}
   Type $\gO(-2,-4)$&
   $(d_1,d_2)=(-2,-4)$ and $(\mu_1^{\#},\mu_2^{\#})=(0,0)$\;, \\
   Type $\gO(-2,-3)$&
   $(d_1,d_2)=(-2,-3)$ and $(\mu_1^{\#},\mu_2^{\#})=(0,1)$ or $(1,0)$\;,\\
   Type $\gO(-1,-4)$&
   $(d_1,d_2)=(-1,-4)$ and $(\mu_1^{\#},\mu_2^{\#})=(1,0)$\;,\\
   Type $\gO(-1,-3)$&
   $(d_1,d_2)=(-1,-3)$ and $(\mu_1^{\#},\mu_2^{\#})=(1,1)$\;,\\
   Type $\gO(-3,-3)$&
   $(d_1,d_2)=(-1,-3)$ and $(\mu_1^{\#},\mu_2^{\#})=(0,0)$\;,\\
   Type $\gO(-2,-2)$&
   $(d_1,d_2)=(-2,-2)$ and $(\mu_1^{\#},\mu_2^{\#})=(1,1)$\;.
\end{tabular}
\end{center}
Since the surface has genus $0$, we can set 
$\overline{M}_0=\C\cup\{\infty\}$ and 
$M=\C\cup\{\infty\}\setminus\{p_1,p_2\}$.  Since $\pi_1(M)$ is commutative, 
all surfaces of these types are reducible.

\subsubsection*{Surfaces of type \boldmath{$\gO(-3,-3)$}}  
There exists a minimal surface in $\R^3$ of class $\gO(-3,-3)$ with total 
absolute curvature $8\pi$ \cite{Lopez}.  The following is proven in 
Section~\ref{sec:deform}:  
\begin{proposition}\label{thm:go33}
  By deforming the minimal surface of type $\gO(-3,-3)$ in $\R^3$, 
  one obtains a one-parameter family of\/ \cmcone{} surfaces of
  type $\gO(-3,-3)$ with dual total absolute curvature $8\pi$.  
\end{proposition}

\subsubsection*{Surfaces of type \boldmath{$\gO(-2,-4)$}}  
In this case, by \eqref{eq:r-r} and \eqref{eq:ord-umb}, such a surface has 
two distinct umbilic points of order $1$.  
We may set the ends to be $(p_1,p_2)=(0,\infty)$ 
and the umbilic points to be $(q_1,q_2)=(1,q)$, $q \in \C \setminus \{0,1\}$, 
on $\C\cup\{\infty\}$.  Then we may assume 
\begin{equation}\label{eq:go24-GQ}
 G=\left( \frac{z-q}{z-1} \right)^2 \; , \qquad 
 Q=\frac{\theta (z-1) (z-q)}{z^2} dz^2 \qquad
   (\theta\in\C\setminus\{0\}) \; .  
\end{equation} 
For such $G$ and $Q$, the roots of the indicial equation of 
\eqref{eq:ode} at $z=0$ are 
\[
    \lambda_1 =
    \frac{1}{2}\left(1+\sqrt{1-4\theta q}\right),\qquad
    \lambda_2 =
    \frac{1}{2}\left(1-\sqrt{1-4\theta q}\right).
\]
Then, by \ref{item:prop:exist:2} of Proposition~\ref{prop:exist}, we have
\begin{theorem}\label{thm:go24}
Let $\theta\in\C\setminus\{0\}$ and $q\in\C\setminus\{0,1\}$ be complex
numbers such that 
\[
   \sqrt{1-4\theta q}\in\R\setminus\Z.
\]
Then there exists a conformal $\Hyp^1$-reducible \cmcone{} immersion 
$f\colon{} \C \setminus\{0\}\to H^3$ of type 
$\gO(-2,-4)$ with $\TA(f^{\#})=8\pi$ whose hyperbolic Gauss map 
and Hopf differential are as in \eqref{eq:go24-GQ}.
Moreover, all $\Hyp^1$-reducible surfaces with $\TA(f^{\#})=8\pi$ of type 
$\gO(-2,-4)$ are given in this manner.
\end{theorem}

It only remains to consider the $\Hyp^3$-reducible case:  

\begin{theorem}\label{thm:go24b}
 Let $s \in \R$ such that $\sqrt{1-4 s} \geq 2$ is an integer.  
 Then there exists at least $1$ and at most $\sqrt{1-4 s}$ conformal 
 $\Hyp^3$-reducible \cmcone{} immersions 
 $f\colon{}\C \setminus\{0 \} \to H^3$ of type 
 $\gO(-2,-4)$ with $\TA(f^{\#})=8\pi$ whose hyperbolic Gauss map 
 and Hopf differential are as in \eqref{eq:go24-GQ}. 
 Moreover, all $\Hyp^3$-reducible surfaces with $\TA(f^{\#})=8\pi$ of 
 type $\gO(-2,-4)$ are given in this manner.  
\end{theorem}
\begin{proof}
 For $G$ and $Q$ in \eqref{eq:go24-GQ}, equation \ref{eq:E-1} becomes 
\begin{equation}\label{eq:go24-h3}
       z^2 X'' + z \left\{2+\frac{4z}{1-z}\right\} X'+  
       \left\{\theta(z-1)(z-q)\right\}X = 0 \; .
\end{equation}
 By Lemma~\ref{lem:added} and Proposition~\ref{prop:exist}, it is enough
 to show that there exists data $(G,Q)$ such that the difference of
 the roots of the indicial equation of \eqref{eq:go24-h3} at $z=0$ 
 is an integer and the log-term vanishes.

 The coefficients of \eqref{eq:go24-h3} expand as
\[
    z\left\{2+\frac{4z}{1-z}\right\} = z\left\{ 2 + 
          4\sum_{j=1}^{\infty}z^j\right\}      \quad
    \text{and}     \quad
    \theta(z-1)(z-q)=\theta q -\theta(1+q) z +\theta z^2
\]
 for $z$ sufficiently close to $0$.   
 Assume the roots $\lambda_1$, $\lambda_2$ of the indicial equation
 of \eqref{eq:go24-h3} satisfy $\lambda_1-\lambda_2=m\in\Z^+$.
 Then
\begin{equation}\label{eq:go24-q0}
     s:=\theta q = \frac{1-m^2}{4}
     \quad \text{and}\quad
     \lambda_2=-\frac{m+1}{2}\qquad (m\geq 2)\;.
\end{equation}
 Let 
\[ 
   \mu_j=
   \begin{cases}
     \dfrac{1}{j(m-j)} \qquad &(j=1,2,\dots,m-1)\\
     -\dfrac{1}{m}     \qquad &(j=m)
   \end{cases}\;.
\] 
 Then by Proposition~\ref{prop:recursion} in Appendix~\ref{app:log},
 the log-term coefficient $c$ of \eqref{eq:go24-h3} is given by
 $ c = a_m$,
  where $a_0=1$ and
\begin{multline}\label{eq:go24-a}
   a_j = \mu_j 
       \left[
        \left(
         \sum_{k=0}^{j-2}(4k-2m-2)a_k
        \right)
        +\theta a_{j-2}\right.\\
        \left.
        +
        \left(\frac{1}{4}(m+1)(m-9)-4+4j-\theta\right)
        a_{j-1}
       \right]
\end{multline}
for $j=1,\dots m$.
Hence $a_j$  is a polynomial in $\theta$ of order $j$.
We now define $t_0 = 1$, and we define $t_j$ and $u_j$ for 
$j=1,\dots,m$ by 
the relations 
\[
   a_j = t_j \theta^j+u_j \theta^{j-1} +\dots \qquad
  (j=1,2,\dots,m)\;.
\]
It follows that $t_j = -\mu_j t_{j-1}$, and hence $t_m \neq 0$.  
Then, defining $\Lambda_j := u_j / t_j$, we also have 
\[ 
   \Lambda_j = \Lambda_{j-1} - \frac{(m+1)(m-9)}{4}-4j+4+(j-1)(m-j+1) 
\] 
for $j=2,\dots,m$.  
Since $\Lambda_1 = -(m+1)(m-9)/4$, we have 
\[ 
   \Lambda_m = \sum_{j=2}^m 
               \left[
                -\frac{(m+1)(m-9)}{4}-4j+4+(j-1)(m-j+1)
               \right] 
             = \frac{m}{12} (49-m^2) \; . 
\]  
If the only roots of the polynomial 
\[
   c = t_m \theta^m+u_m \theta^{m-1} + \dots 
     = t_m( \theta^m +\Lambda_m \theta^{m-1}+\dots)= 0
\]
with respect to $\theta$  are $0$ and $(1-m^2)/4<0$, then 
it follows that $\Lambda_m$ would be nonnegative.
However, $\Lambda_m < 0$ for all $m \geq 8$, hence this polynomial 
must have some root $\theta \in \C \setminus \{0,(1-m^2)/4\}$, and then 
$q = (1-m^2)/(4 \theta) \in \C \setminus \{0,1\}$.  
For this $\theta$ and $q$, we have $c=0$, and thus we have at least one 
surface for each $m \geq 8$.  Since $c$ is a polynomial of degree $m$ in 
$\theta$, there are at most $m$ roots, and hence at most $m$ surfaces.  

For $m \leq 7$, one can check by explicitly computing the polynomial 
for $c$ that there is always at least one root 
$\theta \in \C \setminus \{0,(1-m^2)/4\}$.  
\end{proof}

\subsubsection*{Surfaces of type \boldmath{$\gO(-2,-3)$} with  
            \boldmath$\mu_1^{\#}=0$}  
Here, by \eqref{eq:r-r}, there exists only one umbilic point of order $1$.
We set the ends to be $(p_1,p_2)=(1,\infty)$ and the umbilic point to be 
$q_1=0$.  We may assume 
\begin{equation}\label{eq:go23-0-GQ}
  G = z^2, \qquad 
  Q = \frac{\theta z\,dz^2}{(z-1)^2} \qquad (\theta\in\C\setminus\{0\})\; . 
\end{equation}
Then the roots of the indicial equation of \eqref{eq:ode} at
$z=1$ are
\[
    \lambda_1 =
    \frac{1}{2}\left(1+\sqrt{1-4\theta}\right),\qquad
    \lambda_2 =
    \frac{1}{2}\left(1-\sqrt{1-4\theta}\right).
\]
Hence, by Proposition~\ref{prop:exist}, we have
\begin{theorem}\label{thm:go23-0}
Let $\theta\in\R$ such that $\sqrt{1-4\theta} \in \R \setminus \Z$.  
Then there exists a conformal $\Hyp^1$-reducible \cmcone{} immersion 
$f\colon{}\C\setminus\{1,\infty\}\to H^3$ of type 
$\gO(-2,-3)$ with $\TA(f^{\#})=8\pi$ whose hyperbolic Gauss map 
and Hopf differential are as in \eqref{eq:go23-0-GQ}.
Moreover, all $\Hyp^1$-reducible surfaces of type $\gO(-2,-3)$
with $(\mu_1^{\#},\mu_2^{\#})=(0,1)$ and $\TA(f^{\#})=8\pi$ are 
given in this manner.
\end{theorem}

Now we will show that there are no $\Hyp^3$-reducible surfaces of this type, 
by showing that the log-term coefficient at $z=1$ of \eqref{eq:E-1}
cannot be zero.  
With the data as in \eqref{eq:go23-0-GQ}, equation \ref{eq:E-1} becomes 
\[ 
    (z-1)^2 X''+2 (z-1) X' + \theta (1+(z-1)) X = 0 \; , 
\] 
and so $p_0=2$, $q_0=q_1=\theta$, $p_j=0$ for $j \geq 1$, and 
$q_j=0$ for $j \geq 2$, where the $p_j$ and $q_j$ are as defined in 
\eqref{eq:coef}. 
Therefore, by Proposition \ref{prop:recursion}, we have 
$c=-\theta^m/(m! (m-1)!) \neq 0$. 

\subsubsection*{Surfaces of type \boldmath{$\gO(-2,-3)$} with 
\boldmath$\mu_1^{\#}=1$}  
In this case, we set the ends to be $(p_1,p_2)=(0,\infty)$ and the only 
umbilic point to be $q_1=1$.  Then we may assume 
\begin{equation}\label{eq:go23-1-GQ}
  G = \left( \frac{z-1}{z} \right)^2 \;, \qquad
  Q = \frac{\theta (z-1)\, dz^2}{z^2} \qquad
  (\theta\in\C\setminus\{0\})\;.
\end{equation}
Thus the roots of the indicial equation of \eqref{eq:ode} at $z=0$ are
\[
    \lambda_1 =
    \frac{1}{2}\left(1+\sqrt{4+4\theta}\right),\qquad
    \lambda_2 =
    \frac{1}{2}\left(1-\sqrt{4+4\theta}\right).
\]
So, by Proposition~\ref{prop:exist}, we have
\begin{theorem}\label{thm:go23-1}
Let $\theta \in \R$ such that $\sqrt{4+4\theta} \in \R \setminus \Z$.  
Then there exists a conformal $\Hyp^1$-reducible \cmcone{} immersion 
$f\colon{} \C \setminus \{ 0 \}\to H^3$ of type 
$\gO(-2,-3)$ with $\TA(f^{\#})=8\pi$ whose hyperbolic Gauss map 
and Hopf differential are as in \eqref{eq:go23-1-GQ}.
Moreover, all $\Hyp^1$-reducible surfaces of type $\gO(-2,-3)$
with $(\mu_1^{\#},\mu_2^{\#})=(1,0)$ and $\TA(f^{\#})=8\pi$ are given in 
this manner.
\end{theorem}

Now we will show that there are no $\Hyp^3$-reducible surfaces of this 
type as well, again by showing that a log-term coefficient cannot be zero.  
With $G$ and $Q$ as in \eqref{eq:go23-1-GQ}, equation \ref{eq:E-2}
 becomes 
\[ 
     z^2 X'' - z X' + \theta (z-1) X = 0 \; , 
\] 
and so $p_0=-1$, $-q_0=q_1=\theta$, $p_j=0$ for $j \geq 1$, and 
$q_j=0$ for $j \geq 2$, where the $p_j$ and $q_j$ are as defined in 
\eqref{eq:coef}. Hence again, by Proposition~\ref{prop:recursion}, 
we have $c \neq 0$. 

\subsubsection*{Surfaces of type  \boldmath{$\gO(-1,-4)$}}
We set the ends to be $(p_1,p_2)=(0,1)$ and the single umbilic point to be 
$q_1=\infty$, then we may assume 
\begin{equation}\label{eq:go14-GQ}
  G = z^2\;, \qquad
  Q = \frac{\theta \,dz^2}{z(z-1)^4}\qquad(\theta\in\C\setminus\{0\})\;.
\end{equation}
The roots of the indicial equation of \eqref{eq:ode} for 
such $G$ and $Q$ at $z=0$ are $3/2$ and $-1/2$.  
Then, by Lemma~\ref{m=2}, the log-term coefficient at $z=0$ 
vanishes if and only if $\theta=-4$.
Thus
\begin{theorem}\label{thm:go14}
 Any complete \cmcone{} immersion of type 
 $\gO(-1,-4)$ with $\TA(f^{\#})$ $=8\pi$ is congruent to an $\Hyp^3$-reducible 
 \cmcone{} immersion $f\colon{}M=\C\cup\{\infty\} \setminus \{0,1\} \to H^3$
 with hyperbolic Gauss map and Hopf differential 
\[
  G = z^2\;, \qquad
  Q = \frac{-4 \,dz^2}{z(z-1)^4}\;.
\]
\end{theorem}

\subsubsection*{Surfaces of type \boldmath{$\gO(-1,-3)$}}
In this case, there are no umbilic points, by \eqref{eq:r-r}.
Then, if we set the ends to be $(p_1,p_2)=(0,\infty)$, we may assume 
\[
   G = z^2\;,\qquad Q = \frac{\theta}{z}\,dz^2\qquad
   (\theta\in\C\setminus\{0\})\;.
\]
The roots of the indicial equation of \eqref{eq:ode}
at $z=0$ are $3/2$ and $-1/2$, and the log-term coefficient vanishes if 
and only if $\theta=0$, by \eqref{m=2}.  So this case is impossible.

\subsubsection*{Surfaces of type \boldmath{$\gO(-2,-2)$}}  
Here again there are no umbilic points, by \eqref{eq:r-r}.  
If we set the ends to be $(p_1,p_2)=(0,\infty)$, we may assume 
\[
    G = z^2\;,\quad Q=\frac{\theta}{z^2}\,dz^2 \qquad
   (\theta\in\C\setminus\{0\})\;.
\]
Then the solution $g$ of the equation $S(g)-S(G)=2Q$ is  
\[
    g = a z ^ {\mu} + b\; ,\qquad a\in\C \setminus \{0\}\;,\quad b\in \C
         \qquad\text{and}\quad\mu=\sqrt{1-4\theta} \; . 
\]
Hence the function $g$ satisfies $\rho_g(\tau) \in \SU(2)$ for all 
$\tau \in \pi_1(\C \setminus \{0\})$ if and only if 
$\mu\in\Z$, or $\mu\in\R$ and $b=0$.  In 
the second case, the surface is a double cover of a catenoid cousin.
The first case is a warped catenoid cousin with $m=2$ in Theorem 6.2 of 
\cite{uy1} (see also \cite{ruy4}).  

\subsection*{The case \boldmath $(\gamma,n)=(0,1)$}
In this case, we can set $M = \C$.
Since $M$ is simply connected, we have no period problem.  
By \eqref{eq:ineq-8pi} and \eqref{eq:degree}, $d_1=-5$ or $-6$.  

For the case of $\gO(-5)$, there is one umbilic point, which we may suppose is 
at $q_1=0$.  By \eqref{eq:ineq-8pi}, we have $\mu_1^{\#}=1$, so we may assume 
\begin{equation}\label{eq:go5-GQ}
   G = z^2\;, \qquad Q = \theta z\,dz^2 \qquad (\theta\in\C\setminus\{0\})\;.
\end{equation}

For the case of $\gO(-6)$, there are two umbilic points of order $1$.
Without loss of generality, we can set them to be $(q_1,q_2) =(0,1)$.  
So, since $\mu_1^{\#}=0$, we may assume 
\begin{equation}\label{eq:go6-GQ}
   G = \left(\frac{z-1}{z}\right)^2\;, \qquad 
   Q = \theta z(z-1)\,dz^2 \qquad (\theta\in\C\setminus\{0\})\;.
\end{equation}
\begin{theorem}\label{thm:go5/6}
A \cmcone{} surface of genus zero with one end such that \/$\TA(f^{\#})=8\pi$
is congruent to an immersion $f\colon{}\C\to H^3$ with hyperbolic Gauss map
and Hopf differential as in \eqref{eq:go5-GQ} or \eqref{eq:go6-GQ}.  
Moreover, such a surface is $\Hyp^3$-reducible.
\end{theorem}
\begin{figure}
\begin{center}
\begin{tabular}{cc}
  \includegraphics[width=1.6in]{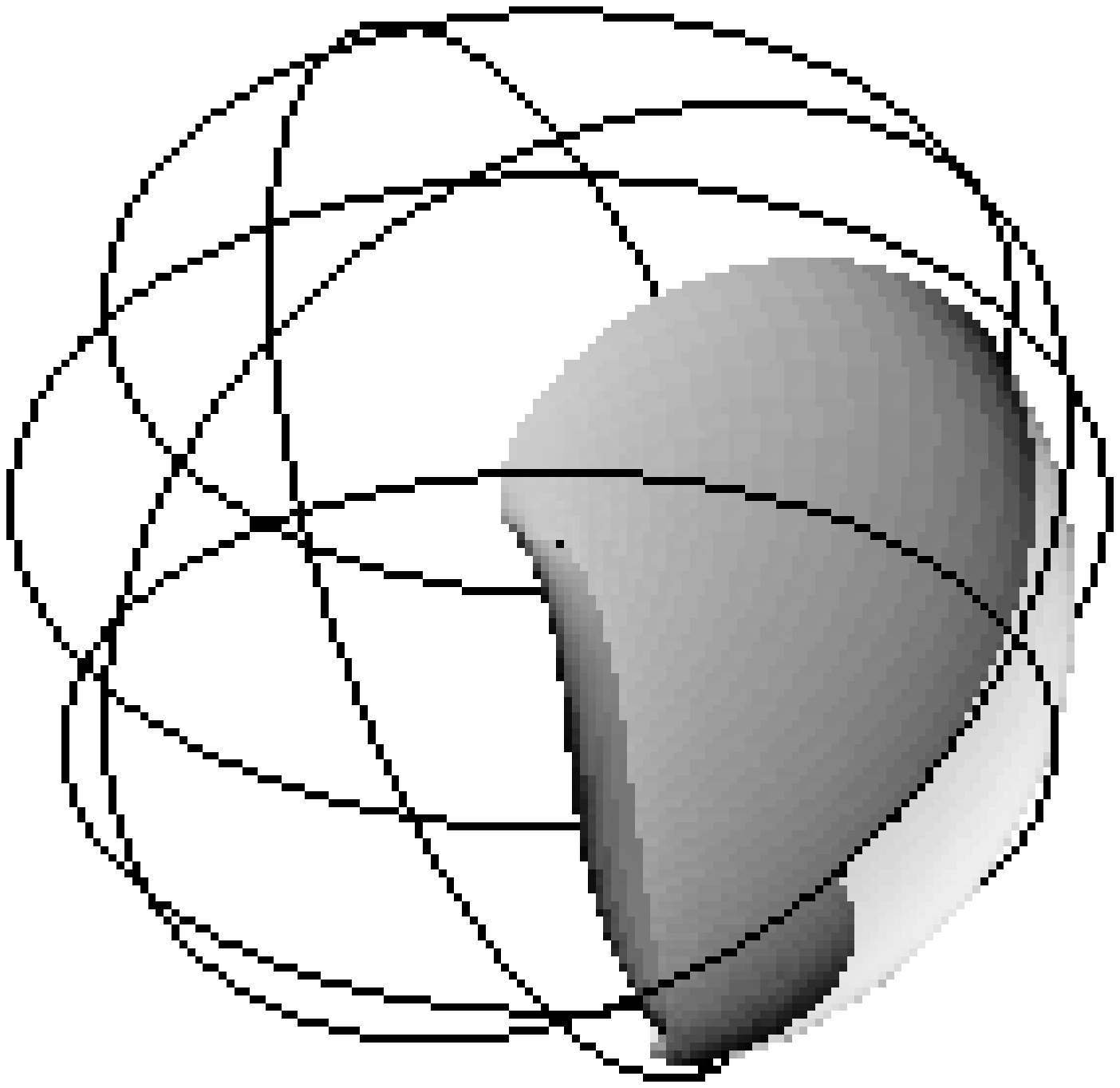} &
  \includegraphics[width=1.5in]{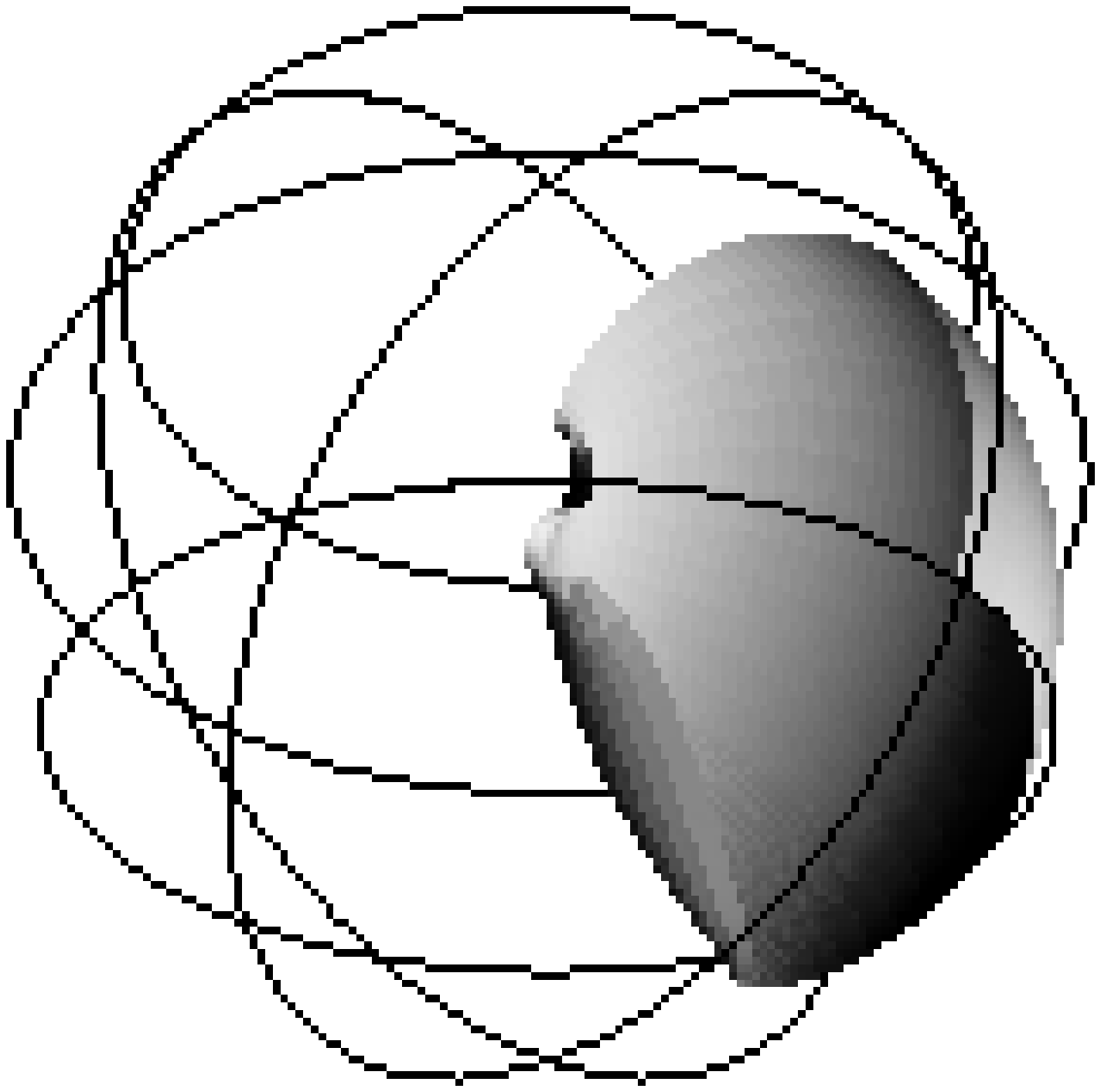}
\end{tabular}
\end{center}
\caption{Genus $0$ and genus $1$ Enneper cousin duals.  
Each surface has a single end that triply wraps around its 
limiting point at the south pole of the sphere at infinity.  
These surfaces are of type $\gO(-4)$ and $\gI(-4)$, and 
have $\TA(f^\#)=4\pi$ and $\TA(f^\#)=8\pi$.  In both cases 
only one of four congruent pieces (bounded by planar geodesics) 
of the surface is shown.}
\end{figure}
\section{Deformation of minimal surfaces to \cmcone{} surfaces}
\label{sec:deform}

In this section, we prove Propositions~\ref{thm:gi4} and \ref{thm:go33}.  
For this, we will need a method from \cite{ruy1} that produces 
a $1$-parameter family of \cmcone{} surfaces in $H^3$ from a corresponding 
minimal surface in $\R^3$, so we describe that method first.  

We start with a complete minimal surface $f_{0}: M\to \R^3$ of 
finite total curvature.  
We require the immersion to be symmetric in the following 
sense, a condition that generically eliminates virtually all minimal 
surfaces, but eliminates none of the better known surfaces, which 
all have symmetries:  
\begin{quote}
{\bf Symmetry condition:}
  There is a disk $D \subset  M$ so that 
  $f_{0}(D)$ is bounded by non-straight planar geodesics.  
\end{quote}
If $f_{0}$ is symmetric with respect to a disk $D$, then $f_{0}(D)$ 
generates the full surface by reflections across planes containing 
the boundary planar geodesics of $\partial f_{0}(D)$, 
by the Schwarz reflection principle \cite{Osserman}.  
Since the surface has finite total curvature, 
it is shown in \cite{ruy1} that the boundary $\partial f_{0}(D)$ is 
contained entirely in only either one plane $P_1$, or in 
two intersecting planes $P_1$, $P_2$, or in three planes 
$P_1$, $P_2$, and $P_3$ in general position.  
Let the boundary planar geodesics of $f_{0}(D)$ contained in $P_j$ 
be called $S_{j,1}$, $S_{j,2}$, \dots, $S_{j,\delta_j}$  
($j=1,\dots,s$, for $s=1$, $2$ or $3$).  

We now define non-degeneracy of the period problems.  
Let $\delta$ be the number of $S_{j,l}$ minus the number of planes 
($\delta= \delta_1+\delta_2+\delta_3-3$ if $s=3$, 
$\delta = \delta_1+\delta_2-2$ if $s=2$, and 
$\delta = \delta_1-1$ if $s=1$).  

\begin{quote} {\bf Nondegeneracy condition:} 
There exists a continuous $\delta$-parameter family of minimal 
disks $f_{0,\nu}(D)$  
(where $\nu$ lies in a small neighborhood
of the origin $\vec{0} \in \R^{\delta}$) such that
  \begin{enumerate}
    \item   $f_{0,\vec{0}}(D) = f_{0}(D)$.
    \item   $\partial f_{0,\nu}(D)= 
            \cup_{j=1}^{s}(\cup_{l=1}^{\delta_j}
              S_{j,l}(\nu))$ holds,  
            and each $S_{j,l}(\nu)$ is 
            a planar geodesic lying in a plane
            $P_{j,l}(\nu)$ parallel to $P_j$.  
    \item   Letting $\Per_{j,l}(\nu)$ 
            ($j=1,\dots,s$, $l=2,\dots,\delta_j$) be the oriented distance 
            between the plane $P_{j,l}(\nu)$ and $P_{j,1}(\nu)$, 
            the map from $\nu$ in $\R^{\delta}$ to 
            $(\Per_{j,l}(\nu))$ in $\R^{\delta}$ is 
            an open map onto a small 
             neighborhood of $\vec{0} \in \R^{\delta}$.  
  \end{enumerate}
\end{quote}

\begin{theorem}[\cite{ruy1}]\label{mainthm} 
If the minimal immersion $f_{0}$ is symmetric and 
nondegenerate, then there exists a one-parameter family of \cmcone{}
surfaces in $H^3$, each of 
whose hyperbolic Gauss map and Hopf differential
coincide with the Gauss map and Hopf differential of
$f_{0,\nu}(D)$ for some $\nu \in \R^{\delta}$.
\end{theorem}

We now consider two applications of this theorem: 
\subsection*{Existence of surfaces of type \boldmath{$\gI(-4)$}
   with \boldmath$\TA(f^{\#}) = 8 \pi$.}

We construct a deformation of 
the Chen-Gackstatter minimal surface defined on the elliptic curves 
\[ 
  \overline{M}_1(\nu_1) = \left\{(z,w) \in (\C \cup \{ \infty \})^2 \; 
  | \; w^2 = z (z-1) (z+\nu_1) \right\} \qquad (\nu_1\in \R^+)\;,
\] 
with the point $p_1$ corresponding to $z=\infty$ removed ($p_1$ will be 
the single end of the surfaces).  Let 
\[
 g= \frac{\nu_2 w}{z} \; , \qquad \omega = \frac{z dz}{w} 
 \qquad (\nu_2\in \R^+)\;. 
 \]  
We choose the fundamental pieces of the surfaces to be the images under 
the Weierstrass representation 
\begin{equation}\label{eq:weier}
  \Re \int_{z_0=0}^z \left( (1-g^2,i(1+g^2),2g \right) \omega
\end{equation}
of the half sheets 
\begin{multline*}
   \left\{ (z,w_1 w_2 w_3) \in 
      \overline{M}_1(\nu_1) \, | \, z \in \C, \Im(z) \geq 0, w_1^2=z, 
        \right. \\
       \left. w_2^2=z-1,w_3^2=z+\nu_1, \arg(w_j) \in [0,\pi), 
     j=1,2,3 \right\} \; . 
\end{multline*}
The fundamental pieces are bounded by 
four planar geodesics, two of which lie in planes parallel to the 
$x_1x_3$-plane and two of which lie in planes parallel to the 
$x_2x_3$-plane.  Thus $\delta=2$.
Note that the period problem is solved, and the Chen-Gackstatter surface 
is produced, if $\nu_1=1$ and $\nu_2=\sqrt{B}$,
where
\[
   B:=\left.\left( \int_0^1 \frac{xdx}{\sqrt{x(1-x^2)}}\right)\right/
         \left( \int_0^1 \frac{(1-x^2)dx}{\sqrt{x(1-x^2)}}\right).
\]
The oriented distance functions (between the two pairs of 
parallel planes containing 
boundary curves of the fundamental pieces) are given by
\begin{align*}
  \Per_1(\nu_1,\nu_2) 
       &= \int_0^1 \left( 1-\nu_2^2 x^{-1} (1-x) (x+\nu_1)
                   \right) 
           \frac{\sqrt{x}\, dx}{\sqrt{(1-x)(x+\nu_1)}} \; , \\
  \Per_2(\nu_1,\nu_2) 
       &=
         \int_0^1 \left( 1-\nu_1 \nu_2^2 x^{-1} (1-x) (x+
         \frac{1 }{ \nu_1}) \right) \frac{\sqrt{\nu_1} \sqrt{x}\, dx }
           { \sqrt{(1-x)(x+(1/\nu_1))}} \; . 
\end{align*}
To see that the period problem is nondegenerate, 
it is sufficient to check that the Jacobian matrix 
$(\partial({\Per}_1,{\Per}_2)/\partial(\nu_1,\nu_2))$ has nonzero 
determinant at $(\nu_1,\nu_2)=(1,\sqrt{B})$.  It is easy to check 
that $|\partial {\Per}_1/\partial \nu_2| = 
|\partial {\Per}_2/\partial \nu_2| \neq 0$ at $(\nu_1,\nu_2)= 
(1,\sqrt{B})$.  Since 
\begin{align*} 
  \left. \frac{\partial {\Per}_1 }{ \partial \nu_1} 
      \right|_{(\nu_1,\nu_2)  = (1,\sqrt{B})} &= 
  \int_0^1 \frac{x+B (1-x^2) }{ 2 (x-1) (1+x)^2 \sqrt{x}} \sqrt{1-x^2}\, dx 
\; , \\  \left. \frac{\partial \Per_2 }{ \partial 
       \nu_1} \right|_{(\nu_1,\nu_2) = (1,\sqrt{B})} &= 
  \int_0^1 \frac{-x(x+2)+B (2+3x)(1-x^2) }{ 2 (x-1) (1+x)^2 
  \sqrt{x}} \sqrt{1-x^2}\, dx \; , 
\end{align*} 
we have 
\[ 
     \left| \frac{\partial \Per_1 }{ \partial 
       \nu_1} \right| \neq \left| \frac{\partial \Per_2 }{ \partial 
           \nu_1} \right| 
\]
at $(\nu_1,\nu_2)=(1,\sqrt{B})$.  
Thus the determinant of the 
Jacobian is nonzero, and the period problem is nondegenerate.  Hence 
Theorem \ref{mainthm} implies existence of associated \cmcone{} surfaces 
in $H^3$ of type $\gO(-4)$.  Furthermore, as Theorem \ref{mainthm} also 
implies that the hyperbolic Gauss maps will be $\nu_2 w/z$, these 
surfaces have dual total absolute curvature $8\pi$.  

\subsection*{Existence of surfaces of type \boldmath{$\gO(-3,-3)$}
with \boldmath$\TA(f^{\#}) = 8 \pi$.}

Let $M=\C\cup\{\infty\}\setminus\{0,\infty\}$ and
\begin{equation}\label{eq:data-O(-3,-3)}
    g=\frac{2 z^2+2 az-a^2-1}{2 (z+1)}+\nu \; ,
    \quad \text{and} \quad
    \omega = \frac{(z+1)^2}{z^3}\,dz \; ,
\end{equation}
where $a,\nu \in \R$.
  
When $\nu=0$, the Weierstrass representation \eqref{eq:weier}
determines a minimal immersion $f_{0}\colon{}M\to\R^3$ with finite 
total curvature of type $\gO(-3,-3)$ (\cite[Theorem~4]{Lopez}).  For the 
metric to be nondegenerate at $z=-1$, we must assume $a \neq -1 \pm 
\sqrt{2}$.  

Since the Hopf differential $Q=\omega\,dg$ satisfies 
$\overline{Q(\bar z)}=Q(z)$, these minimal surfaces each 
have two planar geodesics that are the images of the positive and negative 
real axes of $\C$ under the Weierstrass representation \eqref{eq:weier}, 
and their fundamental pieces are the images of the upper half
plane of $\C$ under \eqref{eq:weier}.  The two planar geodesics comprise 
the boundaries of each of the fundamental pieces, and both lie in planes 
parallel to the $x_1x_3$-plane, since $g$ is real-valued on the real axis.
So $\delta=1$, and the oriented distance between the two planes containing 
the two geodesics is 
\[
   \Per(\nu):=\Re \left( 2\pi i \Res_{z=0} i(1+g^2) \omega \right) 
          =-2 \pi \nu (2  +2  a + \nu) \; , 
\]
so $d \Per(\nu)/d\nu$ is nonvanishing 
at $\nu=0$ when $a \neq -1$.  
Thus Theorem \ref{mainthm} implies existence of a $1$-parameter family of 
\cmcone{} surfaces of type $\gO(-3,-3)$ in $H^3$ for each $a \neq 
-1,-1 \pm \sqrt{2}$ with dual total absolute curvature $8\pi$ (as $g$ has 
degree $2$).  

\appendix
\begingroup
\renewcommand{\thesection}{\Alph{section}}
\section{}
\label{app:log}
\endgroup
\small
Here we review some elementary facts in the theory of
linear ordinary differential equations.  Define a differential operator 
\begin{equation}\label{eq:op}
   L[u]:=z^2 u''+ z p(z)u'+ q(z)u \qquad \left('=\frac{d}{dz}\right)\;.
\end{equation}
In this note, we shall consider the solution of the ordinary
differential equation with a regular singularity at the origin:
\begin{equation}
\label{eq:lineq}
   L[u]= 0 \; , 
\end{equation}
where
\begin{equation}
\label{eq:coef}
  p(z)=\sum_{j=0}^{\infty}p_jz^j\;,\qquad  
  q(z)=\sum_{j=0}^{\infty}q_jz^j\;.
\end{equation}
It is well-known (and we will see it in this appendix) 
that \eqref{eq:lineq} has 
two linearly independent solutions $\{X_1,X_2\}$ of the form 
  \[
     X_1 = z^{\lambda_1} \sum_{j=0}^{\infty}\eta_{1,j} z^j\;,\qquad
     X_2 = \left( z^{\lambda_2} \sum_{j=0}^{\infty}\eta_{2,j} z^j \right) + 
           c \,X_1 \log z\;,
  \]
where $\eta_{1,0} \neq 0$ and $\eta_{2,0} \neq 0$, and where 
$\lambda_1$ and $\lambda_2$ are given by
\begin{equation}
\label{eq:exp}
  \lambda_1 = \frac{1}{2}\left\{(1-p_0)+m\right\}\;,\quad
  \lambda_2 = \frac{1}{2}\left\{(1-p_0)-m\right\}\;,\quad
  m=\sqrt{(1-p_0)^2-4q_0}\;.
\end{equation}
The coefficient $c$ is called {\it the log-term coefficient}
of differential equation \eqref{eq:lineq}, which
may be nonzero only when $\lambda_1-\lambda_2\in \Z$.

\smallskip
We shall give a method for computing the coefficient $c$.
First, we shall describe
two linearly independent solutions $X_1$, $X_2$ as a formal power series.  
If we find a solution of \eqref{eq:lineq} as a formal power series, 
a well-known existence theorem from the theory of ordinary differential 
equations says that it will converge in a sufficiently small 
neighborhood of the origin \cite{cl}.  So the formal treatment is sufficient
for the computation of $c$.  

For a complex variable $\lambda$, define rational functions
$\zeta_j(\lambda)$ for non-negative integers $j$ as
\begin{equation}\label{eq:zeta}
   \zeta_0(\lambda)=1, \quad \text{and}\quad
   \zeta_{j}(\lambda)=  -\frac{1}{\varphi(\lambda+j)}
       \sum_{k=0}^{j-1}r_{j,k}(\lambda)\zeta_k(\lambda)\quad
           (j=1,2,\dots)\;,
\end{equation}
where
\[
    \varphi(t)=t(t-1)+tp_0+q_0,\qquad
    r_{j,k}(\lambda)=(\lambda+k)p_{j-k}+q_{j-k}\;,
\]
and we set 
\begin{equation}\label{eq:formal-X}
  X(\lambda):=z^{\lambda}\sum_{n=0}^{\infty}\zeta_n(\lambda)z^n.
\end{equation}
  Applying the operator $L$ to $X(\lambda)$, we have
\begin{equation}\label{eq:operated}
    L[X(\lambda)]=z^{\lambda}\left\{
                  \varphi(\lambda)+\sum_{j=1}^{\infty}
                    \left(
                      \varphi(\lambda+j)\zeta_j(\lambda)+
                      \sum_{k=0}^{j-1} r_{j,k}(\lambda)\zeta_k(\lambda)
                    \right)z^j
                 \right\}=z^{\lambda}\varphi(\lambda)
\end{equation}

The quadratic equation 
\begin{equation}\label{eq:indicial}
    \varphi(t)=t(t-1)+tp_0+q_0=0
\end{equation}
is called the {\it indicial equation\/} of the equation
\eqref{eq:lineq}, and we denote the solutions of \eqref{eq:indicial}
by $\lambda_1$ and $\lambda_2$.

First, we consider the case $\lambda_1-\lambda_2\not\in\Z$.
In this case, $\varphi(\lambda_l+j)\neq 0$ ($l=1,2$) for any positive
integer $j$, and then $\zeta_j(\lambda_l)$ ($l=1,2$) in \eqref{eq:zeta}
are all well-defined.
Moreover, by \eqref{eq:operated}, $X_1:=X(\lambda_1)$ and 
$X(\lambda_2)$ are linearly independent solutions of \eqref{eq:lineq}.

Next, assume $m:=\lambda_1-\lambda_2$ is a non-negative integer.
Since $\varphi(\lambda_1+j)\neq 0$ for any positive integer $j$, 
$X_1:=X(\lambda_1)$ is a well-defined power series and a solution of
\eqref{eq:lineq}.
\subsubsection*{The case \boldmath{$m=0$}}
   Assume $\lambda_1=\lambda_2$.  
   Since $\varphi(\lambda_1+j)\neq 0$ for any positive integer $j$, 
   $\lambda=\lambda_1$ is not a pole of $\zeta_j(\lambda)$ for each $j$.
   Hence
\[
     \zeta_j(\lambda_1)\qquad\text{and}\qquad
     \left.\frac{\partial}{\partial 
     \lambda}\right|_{\lambda=\lambda_1}\zeta_j(\lambda)
     \qquad (j=0,1,2,\dots)
\]
   are well-defined.
   Let
\begin{equation}\label{eq:sol-m0}
    X_2:=\left.\frac{\partial}{\partial\lambda}\right|_{\lambda=\lambda_1}
               X(\lambda)
        = z^{\lambda_1}
           \sum_{n=0}^{\infty}\left(
           \left.\frac{\partial}{\partial\lambda}\right|_{\lambda=\lambda_1}
           \zeta_n(\lambda)
	   \right) z^n + 
           X_1 \cdot \log z \; . 
\end{equation}
\begin{proposition}\label{prop:m=0}
  If $m=\lambda_1-\lambda_2=0$, $X_2$ in \eqref{eq:sol-m0} is a solution
  of \eqref{eq:lineq}.
  Moreover, the log-term coefficient of \eqref{eq:lineq} never vanishes.
\end{proposition}  
\begin{proof} 
  It is enough to show that $X_2$ is a solution of \eqref{eq:lineq}.
  In fact, by \eqref{eq:operated}, 
\[
    L[X_2] = \left.\frac{\partial}{\partial\lambda}\right|_{\lambda=\lambda_1}
             L[X(\lambda)]
           = z^{\lambda_1}
             \left.\frac{\partial}{\partial\lambda}\right|_{\lambda=\lambda_1}
               \varphi(\lambda)
             +z^{\lambda_1}\varphi(\lambda_1)\log z=0 \; , 
\]
  because $\varphi(\lambda)=(\lambda-\lambda_1)^2$.
\end{proof}
\subsubsection*{The case \boldmath{$m>0$}}
  Assume $m=\lambda_1-\lambda_2$ is a positive integer.
  Since $\varphi(t)=(t-\lambda_2-m)(t-\lambda_2)$, 
  $\varphi(\lambda_2+j)$ does not vanish for each positive integer 
  $j$, except for $j=m$.
  Then $\zeta_j(\lambda)$ has no pole at $\lambda=\lambda_2$ for
  $j=1,2,\dots,m-1$, and 
  may have a pole of order one at $\lambda=\lambda_2$ for $j\geq m$.
  Hence
\[
    \lim_{\lambda\to\lambda_2}\left\{(\lambda-\lambda_2)\zeta_j(\lambda)\right\}
     \qquad
    \text{and}\qquad
    \left.\frac{\partial}{\partial \lambda}\right|_{\lambda=\lambda_2}
    \left[(\lambda-\lambda_2)\zeta_j(\lambda)\right]
\]
  are well-defined.
  Moreover,
\begin{equation}\label{eq:zeta-0}
    \lim_{\lambda\to\lambda_2}
     \left\{(\lambda-\lambda_2)\zeta_j(\lambda)\right\}=0 \qquad
       (j=1,2,\dots, m-1)
\end{equation}
  holds.
  Let
\[
    \xi_j:=\lim_{\lambda\to\lambda_2}
      \left\{(\lambda-\lambda_2)\zeta_{m+j}(\lambda)\right\}\qquad
       (j=0,1,2\dots)
\]
  and set $c:=\xi_0=\lim_{\lambda\to\lambda_2}
  \{(\lambda-\lambda_2)\zeta_m(\lambda) \}$.
  Then by \eqref{eq:zeta} and \eqref{eq:zeta-0},
  we have 
\[
   \xi_0=c\qquad\text{and}\qquad
   \xi_j = \frac{-1}{\varphi(\lambda_2+m+j)}
           \sum_{k=0}^{j-1}r_{j,k}(\lambda_2+m)\xi_k\quad(j=1,2,\dots) \; .
\]
  Comparing this with \eqref{eq:zeta},
  we have $\xi_j = c \zeta_j(\lambda_1)$ ($j=1,2,\dots$), 
  because $\lambda_1=\lambda_2+m$.

  Let 
\begin{equation}\label{eq:sol-m}
    X_2:=\left.\frac{\partial}{\partial\lambda}\right|_{\lambda=\lambda_2}
           \left[\mathstrut(\lambda-\lambda_2)X(\lambda)\right]\;.
\end{equation}
  Then by \eqref{eq:zeta-0}, we have
\begin{align*}
    X_2 
     &=   z^{\lambda_2}
             \left(\sum_{j=0}^{\infty}\xi_j z^{j+m}\right) \log z
           + z^{\lambda_2} \sum_{j=0}^{\infty}
               \left.\frac{\partial}{\partial\lambda}\right|_{\lambda=\lambda_2}
                   \left\{(\lambda-\lambda_2)\zeta_j(\lambda)\right\}z^j\\
     &= c \log z X_1
           + z^{\lambda_2} \sum_{j=0}^{\infty}
               \left.\frac{\partial}{\partial\lambda}\right|_{\lambda=\lambda_2}
                   \left\{(\lambda-\lambda_2)\zeta_j(\lambda)\right\}z^j\;.
\end{align*}
\begin{proposition}\label{prop:m}
  If $m=\lambda_1-\lambda_2$, is a positive integer, 
  $X_2$ in \eqref{eq:sol-m} is a solution of \eqref{eq:lineq}.
  Moreover, the log-term coefficient $c$ of \eqref{eq:lineq} is given by
\begin{equation}\label{eq:log-m}
    c:=\xi_0=\lim_{\lambda\to\lambda_2}
      \left\{(\lambda-\lambda_2)\zeta_m(\lambda)\right\}\;.
\end{equation}
\end{proposition}  
\begin{proof} 
  By \eqref{eq:operated}, 
\[
    L[X_2] = \lim_{\lambda\to\lambda_2} 
      \frac{\partial}{\partial\lambda} \left( 
              z^\lambda(\lambda-\lambda_2)\varphi(\lambda) \right)=0 \; , 
\]
  because $\varphi(\lambda)=(\lambda-\lambda_1)(\lambda-\lambda_2)$.
\end{proof}

We have established the following recursive formula for $c$, 
which follows immediately from equation \eqref{eq:log-m}:  

\begin{proposition}\label{prop:recursion}
  If the difference $m$ of the roots of the indicial equation of 
  \eqref{eq:lineq} is a positive integer, then
  the log-term coefficient $c$ is 
\begin{equation} 
 \label{eq:logcoef3}
    c        = -\frac1{m}\displaystyle\sum_{k=0}^{m-1}
                \left( (\lambda_2+k)p_{m-k}+q_{m-k} \right)a_k\;,
 \end{equation}
where $a_0=1$ and 
\[
    a_j = \frac{1}{j(m-j)}\sum_{k=0}^{j-1}
                \left( (\lambda_2+k)p_{j-k}+q_{j-k} \right)a_k
                \qquad (j=1,2,\dots,m-1) \; . 
\]
\end{proposition}
\begin{proof}
 Since $\varphi(t)=(t-\lambda_2)(t-\lambda_2-m)$,
 $\varphi(\lambda_2+j)\neq 0$ for $j=1,\dots,m-1$ and then
 $a_j = \zeta_j(\lambda_2)$ ($j=1,\dots,m-1$) is well-defined.
 Hence, by \eqref{eq:log-m},
\begin{align*}
    c &= \lim_{\lambda\to\lambda_2}
         \left\{(\lambda-\lambda_2)\zeta_m(\lambda)\right\}\\
      &= \lim_{\lambda\to\lambda_2}
         \frac{-(\lambda-\lambda_2)}{%
               (\lambda+m-\lambda_2)(\lambda-\lambda_2)}
         \sum_{k=0}^{m-1} r_{m,k}(\lambda) \zeta_k(\lambda)\\
      &=  -\frac{1}{m}\sum_{k=0}^{m-1} \left( 
               (\lambda_2+k)p_{m-k}+q_{m-k} \right) a_k \; .
\end{align*}
 This completes the proof.
\end{proof}

Thus, in the case that $p(z)=0$ and $m=1$, $2$, or $3$, 
the solutions of $z^2 u''(z)+q(z) u(z)=0$
have no log-term if and only if
\begin{align}
  q_1                        &=0 \qquad (m=1) \; , \label{m=1}\\
  q_2+(q_1)^2                &=0 \qquad (m=2) \; , \label{m=2}\\
  q_3+q_1q_2+\frac14(q_1)^3  &=0 \qquad (m=3) \; , \label{m=3}
\end{align}
where $q(z)=\sum_{j=0}^\infty q_j z^j$, as in \eqref{eq:coef}.  



\begin{thebibliography}{20}
\bibitem{Bryant}
  R.~Bryant,
  {\itshape Surfaces of mean curvature one in hyperbolic space},
  Ast\'erisque {\bfseries 154--155} (1987), 321--347.
\bibitem{cg}
  C.~C.~Chen, F.~Gackstatter,
  {\itshape Elliptische und hyperelliptische Funktionen und
    vollst\"andige Minimalfl\"achen vom Enneperschen Typ},
    Math.~Ann. {\bfseries 259} (1982), 359--369.\
\bibitem{cl}
  E.~A.~Coddington, N.~Levinson,
  {\sc Theory of Ordinary Differential Equations}, McGraw-Hill, 1955.  
\bibitem{chr}
  P.~Collin, L.~Hauswirth and H.~Rosenberg,
  {\itshape The geometry of finite topology surfaces properly embedded 
            in hyperbolic space with constant mean curvature one},
  preprint.
\bibitem{hc}
   A.~Hurwicz and R.~Courant, 
   {\sc Funktionen theorie}, 4.~Auflage,  Springer, 1964.
\bibitem{Lopez}
   F.~J.~Lopez,
   {\itshape The classification of complete minimal surfaces with
        total curvature greater than $-12\pi$},
   Trans.~Amer.~Math.~Soc. {\bfseries 334} (1992), 49--74.

\bibitem{mu}
   C.~McCune and M.~Umehara,
   {\itshape 
   An analogue of the UP-iteration for
constant mean curvature one surfaces in Hyperbolic $3$-space},
   preprint.

\bibitem{Osserman}
    R.~Osserman,
  {\sc A Survey of Minimal Surfaces}, {2nd ed.}, Dover, 1986.
\bibitem{rs}
  W.~Rossman, K.~Sato,
  {\itshape Constant mean curvature surfaces with two ends in hyperbolic
  space}, Experimental Math., {\bfseries 7(2)} (1998), 101--119.
\bibitem{ruy1}
  W.~Rossman, M.~Umehara and K.~Yamada,
  {\itshape Irreducible constant mean curvature $1$ surfaces in
    hyperbolic space with positive genus}, 
  T\^ohoku Math.~J. {\bfseries 49}  (1997), 449--484.
\bibitem{ruy2}
  \bysame,
  {\itshape A new flux for mean curvature $1$ surfaces
    in hyperbolic $3$-space, and applications},
    Proc.~Amer.~Math.~Soc. {\bfseries 127} (1999), 2147--2154.
\bibitem{ruy3}
  \bysame,
  {\itshape Mean curvature $1$ surfaces with low total curvature 
  in hyperbolic $3$-space} (an announcement), 
   to appear in J.A.M.I. proceedings, Advanced Studies in Pure Mathematics.
\bibitem{ruy4}
  \bysame,
  {\itshape Mean curvature $1$ surfaces in hyperbolic
            $3$-space with low total curvature II}, 
  in preparation.
\bibitem{s}
  R. Schoen,
  {\itshape Uniqueness, symmetry, and embeddedness of minimal surfaces},
    J. Diff. Geom. {\bfseries 18} (1982), 791--809.
\bibitem{uy1}
  M.~Umehara and K.~Yamada,
  {\itshape Complete surfaces of constant mean curvature-$1$
       in the hyperbolic $3$-space},
  {Ann. of Math. {\bfseries 137} (1993), 611--638.}
\bibitem{uy2}
  \bysame,
  {\itshape A parametrization of Weierstrass formulae and 
     perturbation of some complete minimal surfaces of
        $\R^3$ into the hyperbolic $3$-space},
   J. reine u.~angew.~Math. {\bfseries 432} (1992), 93--116.
\bibitem{uy3}
   \bysame,
   {\itshape Surfaces of constant mean curvature-$c$
        in $H^3(-c^2)$ with prescribed hyperbolic Gauss map},
   Math. Ann. {\bfseries 304} (1996), 203--224.
\bibitem{uy4}
   \bysame,
   {\itshape Another construction of a \cmcone{} surface in $H^3$},
    Kyungpook Math. J. {\bfseries 35} (1996), 831--849.
\bibitem{uy5}
   \bysame,
   {\itshape A duality on \cmcone{} surface in the hyperbolic $3$-space
        and a hyperbolic analogue of the Osserman Inequality},
    Tsukuba J. Math. {\bfseries 21} (1997), 229-237.
\bibitem{uy6}
   \bysame,
   {\itshape Metrics of constant curvature $1$ with
        three conical singularities on the $2$-sphere}, 
   Illinois J. Math. {\bfseries 44} (2000), 72--94.
\bibitem{Yoshida}
M.~Yoshida, {\it Fuchsian Differential Equations},
Max-Plank-Institut f\"ur Mathematik, Friedr. Vieweg \& Sohn, Bonn 1987. 
\bibitem{Yu}
    Z.~Yu,
   {\itshape Value distribution of hyperbolic Gauss maps},
  Proc.~Amer.~Math.~Soc. 
  {\bfseries 125} (1997), 2997--3001.
\bibitem{Yu2}
   \bysame,
   {\itshape The inverse surface and the Osserman Inequality},
  Tsukuba J. Math. 
  {\bfseries 22} (1998), 575--588.
\end{thebibliography}
\end{document}